%% file: main_arxiv.tex
\documentclass[letterpaper]{article}
\usepackage{times}
\usepackage{helvet} \usepackage{courier}  

\usepackage[table]{xcolor}
\RequirePackage{setspace,amsmath,amsfonts,amssymb,amsthm}
\usepackage{graphicx}
\usepackage{enumerate}
\usepackage{tikz}
\RequirePackage[colorlinks,citecolor=blue,urlcolor=blue]{hyperref}
\usetikzlibrary{patterns}
\usepackage{pgfplots}
\pgfplotsset{compat=1.14}
\usepackage{array}
\usepackage{booktabs} 
\usepackage{siunitx} 
\usepackage{pgfplotstable} 
\usepackage{diagbox}

\newtheorem{theorem}{Theorem}
\newtheorem{lemma}{Lemma}

\usepackage{times}
\usepackage[T1]{fontenc}
\usepackage[utf8]{inputenc}
\usepackage{bm}
\usepackage{natbib}

\graphicspath{{./Figs/}}

\usepackage[plain,noend]{algorithm2e}

\input{macros}

\input{Figs/tkz/def_colors}

\title{Higher criticism for rare and weak non-proportional hazard deviations in survival analysis}
\author{A. Kipnis \\
School of Computer Science, Reichman University \\
\and B. Galili \\
Department of Computer Science, Technion - Israel Institute of Technology \\
\and Z. Yekhini \\
School of Computer Science, Reichman University \\
Department of Computer Science, Technion - Israel Institute of Technology} 
\begin{document}
\maketitle
… 

\begin{abstract}
\input{abstract_text}
\end{abstract}

\section{Introduction}
\label{sec:intro} 
\input{intro.tex}

\section{Power Analysis under an Exponential Decay with Rare and Weak Departures}
\label{sec:analysis}
\input{analysis.tex}

\section{Empirical Results}
\label{sec:simulations}
\input{simulations.tex}

\section{Discussion \label{sec:discussion}}
\input{discussion.tex}


\section{Acknowledgments}
The authors would like to thank Malka Gorfine for useful comments on an earlier version of this manuscript. The work of Alon Kipnis is funded in part by the US-Israel Binational Science Foundation (BSF grant No. 2022124).

\section{Supplementary Material}
The Supplementary Material includes the proofs of Theorems~\ref{thm:HC_powerful}, \ref{thm:HC_powerless}, and \ref{thm:LR_powerless}. The code for the simulations is available at \url{https://github.com/alonkipnis/HCHG}.

\input{proofs.tex}

\bibliography{survival}
\bibliographystyle{agsm}
\end{document}

%% file: macros.tex
\usepackage[normalem]{ulem}

\newcommand{\Pois}{\mathrm{Pois}}

\newcommand{\Bin}{\mathrm{Binomial}}

\newcommand{\Prp}[1]{\Pr\left[#1 \right]}
\newcommand{\reals}{\mathbb R}
\newcommand{\simiid}{\overset{\mathsf{iid}}{\sim}}

\newcommand{\Ncal}{\mathcal{N}}
\newcommand{\Ucal}{\mathcal{U}}
\newcommand{\Qcal}{\mathcal{Q}}

\newcommand{\lognull}{\mathrm{Exp}(2)}

\newcommand{\one}{\mathbf{1}}
\newcommand{\HC}{\mathrm{HC}}
\newcommand{\HCHG}{\mathrm{HCHG}_T}
\newcommand{\LR}{\mathrm{LR}_T}
\newcommand{\HG}{\mathrm{HyG}}
\newcommand{\minP}{\mathrm{minP} }

\newcommand{\Var}{\mathrm{Var}}
\newcommand{\ex}[1]{\ensuremath{\mathbb{E}\left[ #1\right]}}

%% file: Figs/tkz/def_colors.tex
\definecolor{tableau_red}{RGB}{255,87,89}
\definecolor{tableau_blue}{RGB}{91,155,213}
\definecolor{tableau_green}{RGB}{140,198,63}
\definecolor{tableau_orange}{RGB}{255,127,15}
\definecolor{tableau_yellow}{RGB}{255,187,120}
\definecolor{tableau_purple}{RGB}{142,68,173}
\definecolor{tableau_pink}{RGB}{197,27,138}
\definecolor{tableau_brown}{RGB}{140,86,75}
\definecolor{tableau_gray}{RGB}{144,144,144}
\definecolor{tableau_cyan}{RGB}{23,190,207}
\definecolor{lightgray}{RGB}{211,211,211}

\definecolor{HCcolor}{HTML}{d1615d}
\definecolor{LRcolor}{HTML}{5778a4}

%% file: abstract_text.tex
We propose a method to compare survival data based on higher criticism of p-values obtained from many exact hypergeometric tests. The method accommodates non-informative right-censorship and is sensitive to hazard differences in unknown and relatively rare time intervals. It attains much better power against such differences than the log-rank test and its variants. We demonstrate the usefulness of our method in detecting rare and weak non-proportional hazard differences compared to existing tests, using simulations and actual gene expression data. Additionally, we analyze the asymptotic power of our method and other tests under a theoretical framework describing two groups experiencing failure rates that are usually identical over time, except in a few unknown instances where one group's failure rate is higher. Our test's power experiences a phase transition across the plane of rarity and intensity parameters that mirrors the phase transition of higher criticism in two-sample rare and weak normal and Poisson means settings. The region of the plane in which our method has asymptotically full power is larger than the corresponding region for the log-rank test.

%% file: intro.tex
\subsection{Survival data with rare hazard deviations}
\begin{figure}[ht!]
\begin{center}
\includegraphics[scale=0.45, trim=0cm 0cm 0cm 1.4cm, clip=true]{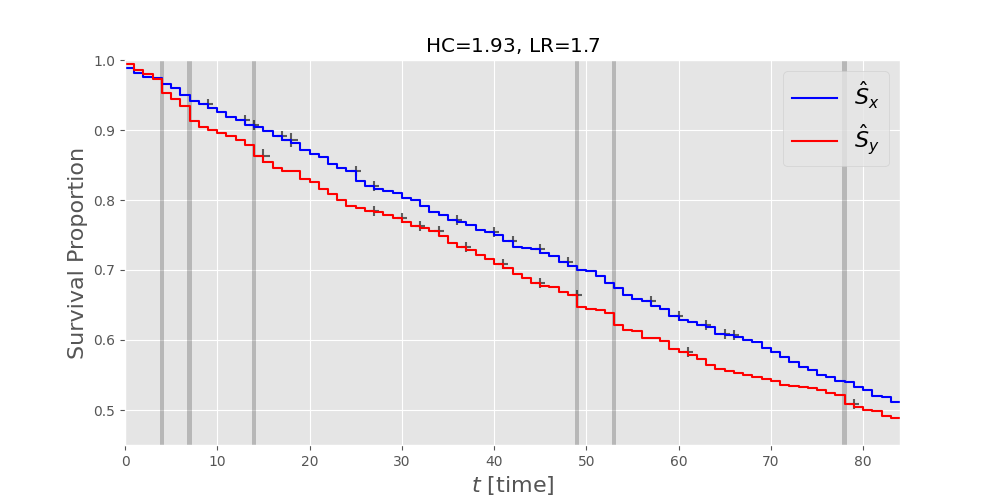}
\include{table_data_intro}
\end{center}
    \caption{Survival data 
    thought to experience temporarily rare excessive hazard in Group $y$ compared to Group $x$. Higher criticism of the hypergeometric p-values (HCHG) indicates an excessive hazard in Group $y$, while the log-rank test does not. Top (figure): Kaplan-Meier curves of the data. Gray bars indicate membership in the set $\Delta^\star$ of time instances providing the best evidence for a global rare hazard difference. Intervals with censoring events are decorated with $+$. Bottom (table): At-risk subjects and events occurring in two groups over $T=84$ time intervals, and the corresponding hypergeometric p-values $\{p_t\}_{t=1}^T$ of \eqref{eq:pvals_def}.
    }
    \label{fig:data_intro}
\end{figure}

Suppose that we have survival measurements from two groups, say, the Control Group $x$ and the Treatment Group $y$. We want to determine whether the treatment significantly affects survival in the sense that the global difference between groups' failure rates goes beyond expected fluctuations. In general, this is a topic studied for many decades with plenty of scientific and industrial applications \citep{kiefer1988economic,armitage2008statistical,kalbfleisch2011statistical}. One notable tool to compare survival data is the log-rank test introduced by \cite{mantel1966evaluation}, which can accommodate right-censorship in the data and is asymptotically equivalent to the likelihood ratio test under the Cox proportional hazard risk model \citep{peto1972asymptotically,kalbfleisch2011statistical,galili2021stability}. Many variations of the log-rank test were proposed to address non-proportional hazard situations \citep{gill1980censoring,harrington1982class,pepe1991weighted,yang2010improved,liu2022log,gorfine2020k}. 
Nevertheless, as we explain below and demonstrate in Table~\ref{tbl:no_discoveries}, existing tools are typically ineffective when hazard differences between groups are rare (sparse). Namely, when the signal separating the two groups corresponds to a few time intervals experiencing excessive or reduced failure rates while these intervals are unknown to us. This paper aims to develop a tool that can reliably detect such temporarily rare and weak hazard departures in survival data and can rigorously handle non-informative right-censored data. As a by-product, the tool can also indicate those time intervals suspected of experiencing increased or decreased hazard. We illustrate these points using example survival data in Figure~\ref{fig:data_intro}. This figure shows the Kaplan-Meier curve of example survival data with significant excessive risk in Group $y$ according to our method, but not according to the log-rank. The gray bars in this figure indicate time intervals thought to experience increased hazard.
\begin{table}[h]
\begin{center}
    \input{table_intro_discovery}
\end{center}
    \caption{True discovery proportion of several tests for survival data in $1,000$ independent random experiments at significance level $\alpha=0.05$. In each experiment, we sample from the rare and weak non-proportional hazard model \eqref{eq:model_full1} below and evaluate several test statistics: HCHG is our newly proposed method, the family of KONP tests was proposed in \cite{gorfine2020k}, FH stands for the family of Fleming-Harrington tests \citep{harrington1982class}, the other tests are Tarone-Ware \citep{tarone1977distribution}, 
    Gehan-Wilcoxon \citep{gehan1965generalized}, and Peto-Prentice \citep{peto1972asymptotically,prentice1978linear}; the expected number of non-null intervals is $4$ out of $T=84$ intervals; all tests are two-sided.
    \label{tbl:no_discoveries}
    }
\end{table}

Below are examples where our method may be particularly useful compared to existing ones.

\subsubsection{Discovering age-specific effects} Suppose that some genetic property may cause disease at certain ages of an organism, but we do not know at what ages the effect will occur. Therefore, to decide whether the effect is significant, we look for differences in the rate of occurrence of the disease across all ages. For example, this situation seems relevant to studying life span quantitative trait loci in Drosophila melanogaster \citep{10.1534/genetics.104.038331}. 
    
\subsubsection{Comparing the rate of decay of radioactive materials} The decay rate of some radioactive materials appears to fluctuate in time due to various causes such as ``space weather'' \citep{milian2020fluctuations}. To identify specific causes, we may compare the survival curve of the material in the exposed environment to that of the same material in a controlled environment. Due to the potential burstiness of space weather, if any effect exists, it may manifest through rare fluctuations in the decay rate of the exposed material. Our test is designed to detect such fluctuations. 

\subsubsection{Identifying possibly opposing temporarily localized hazard trends} As we explain below, our method is also naturally suited to detect situations where hazard differences between the two groups are positive in some time intervals, negative in others, and potentially zero in most. For example, the study \cite{johansson2015family} suggests that the effect of pregnancy on breast cancer risk is positive in the short term and negative in the long term. Since the opposing trends may cancel each other, it is challenging to detect a global hazard difference using methods based on averaging, like the log-rank test and its generalizations. Our method is particularly useful compared to other methods in these situations when the direction and the location of the differences are unknown in advance, so in the example above the concepts ``short" and ``long" can be objectively determined from the data. The situation described here appears to occur also in the studies \cite{tsodikov2002semi,sasaki2005temporal,dekker2008survival,daniels2017examining,gregson2019nonproportional} and in effects of the type ``what doesn't kill you makes you stronger'' \citep{stenton2022effects,mathew2011autophagy,Dongjuan2017}. 

\subsection{Existing methods and rare effects}
The most popular methods for discovering differences or excessive risk in survival data are based on averaging some observed quantities across all event times as in the log-rank test mentioned earlier and its generalizations
\citep{mantel1966evaluation,peto1972asymptotically,gill1980censoring,harrington1982class,pepe1991weighted,galili2021stability,liu2022log,bardo2023methods}. Therefore, it may be the case that intervals of excessive or reduced hazard exist in the data, but these are so rare that they go undetected by these averaging-based approaches -- even if the effect's direction is the same in all non-null intervals which is the typical situation that we address here. More specifically, weighted versions of the log-rank test can be useful against non-proportional hazard alternatives only when the hazard pattern is known in advance; see the additional discussion and references in \cite{yang2010improved} and \cite{chauvel2014tests}. 
Selecting the weights from the data may improve the power against non-proportional hazard alternatives under certain hazard difference models \citep{yang2010improved,chauvel2014tests}, but such adaptive selection necessarily results in loss of power when this function cannot be estimated reliably as in our case of very rare and weak differences. This limitation is well-understood in the context of the Gaussian sequence model \citep{jin2003detecting}. 
%
Additionally, if the hazard in a certain interval differs between the groups, this difference may still be \emph{weak} in the sense that the global effect is undetectable in a Bonferroni analysis involving the significance of individual time intervals \citep{jin2016rare}. For the same reasons, methods based on the maximum of several standardized tests are also ineffective in general (again, unless the hazard pattern is known beforehand) \citep{breslow1984two,self1991adaptive,fleming1987supremum,fleming2013counting}. In contrast, we propose to combine signals from individual time intervals using higher criticism, which is reputed to be effective in detecting rare and individually weak effects \citep{donoho2015special}.

\subsection{Setting}
We have two series of positive integers $\{n_x(t)\}_{t=0}^T$ and $\{n_y(t)\}_{t=0}^T$ of equal length, describing the number of subjects at risk at times $t=1,\ldots,T$ in groups $x$ and $y$, respectively. We are also given the sequences
$\{o_x(t)\}_{t=1}^T$, and $\{o_y(t)\}_{t=1}^T$, describing the number of events occurring in each group over time, so that
\[
o_x(t) \leq n_x(t-1) - n_x(t), \qquad o_y(t) \leq n_y(t-1) - n_y(t),\qquad t=1,\ldots,T,
\]
with equality only if none of the subjects within the corresponding group were censored between time $t-1$ and $t$. We assume that failure and censoring times are independent within each group as in standard log-rank analysis. 

Denote by $c_{\star}(t)$ the number of censored subjects up to time $t$ in group $\star \in \{x,y\}$. The survival proportion (aka estimated survival probability) at time $t$ is 
\[
\hat{s}_{\star}(t) := \frac{n_{\star}(t)}{n_{\star}(0)- c_{\star}(t)},\qquad t=1,\ldots,T.  
\]
The Kaplan-Meier survival curve associated with the group ${\star}$ is the graph 
\begin{align}
    \label{eq:Kaplan_Meier}
    \left\{ \left(t, \hat{s}_{\star}(t)\right) \right\}_{t=1,\ldots,T}, \qquad{\star} \in \{x,y\}.
\end{align}
This curve describes the proportion of at-risk subjects at time $t$ in Group ${\star}$ with censored subjects removed; see Figure~\ref{fig:data_intro} for an example of the Kaplan-Meier curves of survival data.

\subsection{Method description}
\label{sec:method}
We now describe our statistical test for comparing the survival of Group $x$ and Group $y$. Our test uses the higher criticism of p-values obtained from many exact hypergeometric tests as per the explanation below.

\subsubsection{Hypergeometric p-values and Survival Analysis}
The hypergeometric distribution $\HG(M, N, n)$ has the probability mass function 
\begin{align}
    \label{eq:hyg_pmf}
\Prp{\HG(M, N, n)=k} = \frac{\binom{N}{k}\binom{M-N}{n-k}}{\binom{M}{n}},
\end{align}
describing the probability of observing $k$ type-$A$ items in a random sample of $n$ items without replacement from a population of size $M$, initially containing $N$ type-$A$ items. 

Given an observed value $m \in \mathbb N$, the one-sided P-value of the exact hypergeometric test is
\[
p_{\HG}(m; M, N, n) := \Prp{ \HG(M, N, n) \geq m } = \sum_{k=m}^{n} \frac{\binom{N}{k}\binom{M-N}{n-k}}{\binom{M}{n}}.
\]

Back to survival analysis. For every $t=1,\ldots,T$, we evaluate:
\begin{align}
    \label{eq:pvals_def}
    & p_t := p_{\HG}( m_t; M_t, N_t, n_t),
\end{align}
with
\[
m_t = o_y(t), \qquad 
M_t = n_x(t-1) + n_y(t-1),\quad N_t = n_y(t-1),\quad n_t = o_x(t) + o_y(t).
\]
In words, $p_t$ is a P-value under the model proposing that the number of failure events $o_y(t)$ observed in Group $y$ at time $t$ is obtained by sampling without replacement $o_x(t) + o_y(t)$ subjects from a pool of $n_x(t-1) + n_y(t-1)$ subjects, out of which $n_y(t-1)$ subjects are 'at-risk' in Group $y$ at the beginning of the $t$-th interval. The hypergeometric P-value $p_t$ is small if $o_y(t)$ is much larger than the expected number of such events under this model. The table in Figure~\ref{fig:data_intro} shows an example of survival data and the corresponding hypergeometric p-values.

\subsubsection{Higher Criticism}
In this work, we combine the p-values $p_1,\ldots,p_T$ using the HC statistic \citep{donoho2004higher,donoho2008higher}. Specifically, set
\begin{align*}
    \HC_{i;T}(p_1,\ldots,p_T) := \sqrt{T} \frac{i/T - p_{(i)}}{\sqrt{p_{(i)}(1-p_{(i)})}},\quad i=1,\ldots,T,
\end{align*}
where $p_{(1)}\leq \ldots \leq p_{(T)}$ are the ordered p-values observed in the data. The higher criticism statistic of Hyper Geometric p-values (HCHG) is defined as
\begin{align}
    \HCHG := \HC \left(p_1,\ldots,p_T;\gamma_0\right) := \max_{1\leq i \leq T\gamma_0} \HC_{i;T}(p_1,\ldots,p_T)
    \label{eq:HC_def}.
\end{align}
Here $\gamma_0 \in (0,1]$ is a tunable parameter that does not change the large sample properties of $\HCHG$ under either hypothesis \citep{donoho2004higher}. Our test rejects the null hypothesis of equal population survival rates for large values of $\HCHG$. 

It might be reasonable to replace higher criticism with other statistics that are sensitive to rare effects like the Berk-Jones statistics \citep{moscovich2016exact} or the family of phi-divergence statistics introduced in \cite{jager2007goodness}. We focus on higher criticism mainly due to its simplicity and the inherent thresholding mechanism that identifies intervals suspected of excessive hazard as we discuss later on. 

\subsubsection{Critical test values}
 \label{sec:critical_values}
 We are interested in characterizing critical values for testing using $\HCHG$ at a prescribed significance level $\alpha$. We propose to obtain these values by simulating samples of $\HCHG$ under a proposed null model $H_0$. Namely, given a large simulated sample, we consider the empirical $1-\alpha$ quantile as an estimate of the 
 true quantile 
 \[
q_0^{1-\alpha}(\HCHG) := 
\inf \{q \,:\, \Prp{{\HCHG \leq q|H_0}} \geq 1-\alpha \}.
 \]
 This estimate of $q_0^{1-\alpha}(\HCHG)$ serves as the critical value above which we reject the null at level $\alpha$. 

Our experience shows that a test based on the simulated $q_0^{1-\alpha}(\HCHG)$ has much better power than a test based on simulating the $1-\alpha$ quantile of $\HC$ of p-values that are uniformly distributed over $(0,1)$. Indeed, because the data is discrete, the distribution under the null of the hypergeometric p-values is in many cases significantly stochastically larger than uniform hence the null values of $\HCHG$ can be significantly smaller than those obtained when the p-values follow a uniform distribution. Consequently, the $1-\alpha$ quantile of a sample of HC of uniform p-values is overly conservative for an $\alpha$-level test. The max Brownian bridge distribution to which HC asymptotes also leads to overly conservative critical values due to the same reason and also due to the slow convergence of HC to its asymptotic distribution from below \citep{donoho2004higher,gontscharuk2015intermediates,moscovich2016exact}. 
A standard decision-theory solution to improve power while controlling the level when data is discrete is to randomize tests so that the p-values are uniform under any null model for the data 
\citep[p. 101]{cox1979theoretical}, \citep{habiger2011randomised,chen2020false}. However, in large-scale multiple testing, randomizing individual tests introduces additional issues concerning decision-making and interpreting instances of departure \citep{kulinskaya2009fuzzy,efron2012large}.  
Therefore, in analyzing real data, we use non-randomized hypergeometric tests and simulate the null distribution of $\HCHG$ by randomly assigning group membership to subjects, ignoring the actual group membership associated with biological traits. This practice, known as the permutation approach to test calibration, has been recently discussed in contexts related to the sparse signal setting of this paper
\citep{arias2018distribution,kim2022minimax,dobriban2022consistency,stoepker2024anomaly,stoepker2023sparse}.

\subsection{Effect direction}
In most applications we are interested in testing whether the risk in Group $y$ (say) is larger than that of Group $x$, hence we use one-sided p-values in \eqref{eq:pvals_def}. Replacing the one-sided p-values with two-sided ones may be justified when we are interested in detecting a global effect that can go either way, coherently or incoherently, over time. 
Even with one-sided p-values, a significant value of $\HCHG$ might remain significant after replacing the roles of $x$ and $y$. If this situation occurs, then our data may experience a global effect that changes over time, e.g., excessive hazard in the short term and reduced hazard in the long term. Effects of this type were studied in \cite{dekker2008survival,johansson2015family,gregson2019nonproportional,daniels2017examining,mathew2011autophagy,Dongjuan2017}. 

Additionally, we are often interested in a \emph{strictly one-sided effect} in which one group experiences an excessive hazard compared to the other group. We declare that ``Group $y$ experiences an increased failure rate compared to Group $x$" only if HCHG rejects when testing against increased mortality in Group $y$ but does not reject when testing against increased mortality in Group $x$. If each HCHG test is of significance level $\alpha$, the combined test clearly rejects at significance level $\alpha$ or smaller. 

We summarize the general one-sided procedure in Algorithm~\ref{alg:1} and its strictly one-sided variant in Algorithm~\ref{alg:2}

\begin{algorithm}[ht]
\KwResult{reject/not reject $H_0$}
\KwData{Survival data $\{n_x(t), n_y(t), o_x(t), o_y(t) \}_{t=0}^T$.}
\For{$t \in 1,\ldots,T$}{
    $M_t \leftarrow n_x(t) + n_y(t)$\;
    $n_t \leftarrow o_x(t) + o_y(t)$\;
    $p_t \leftarrow p_{\HG}\left( o_y(t); M_t, n_y(t), n_t\right)$\;
}
$\HCHG \leftarrow \HC \left(p_1,\ldots,p_T \right)$

\eIf{$\HCHG > q_0^{1-\alpha}(\HCHG)$ }
{\text{reject $H_0$}}
{\text{do not reject $H_0$}}
\caption{TestHCHG. Testing against an excessive risk in Group $y$.
\label{alg:1}  }
\end{algorithm}

\begin{algorithm}[ht]
\KwResult{reject/not reject $H_0$}
\KwData{Survival data $\{n_x(t), n_y(t), o_x(t), o_y(t) \}_{t=0}^T$.}
$\mathcal S_{(x,y)} \leftarrow \{n_x(t), n_y(t), o_x(t), o_y(t) \}_{t=0}^T$\;
$\mathcal S_{(y,x)} \leftarrow \{n_y(t), n_x(t), o_y(t), o_x(t) \}_{t=0}^T$\;
\eIf{($\mathrm{TestHCHG}$($\mathcal S_{(x,y)}$) rejects $H_0$) \& ($\mathrm{TestHCHG}$($\mathcal S_{(y,x)}$) does not reject $H_0$)} 
    {\text{reject $H_0$}}
{\text{do not reject $H_0$}}
\caption{Testing against a strictly one-sided effect of excessive risk in Group $y$ 
\label{alg:2}  }
\end{algorithm}







\subsection{Identifying instances of departure \label{sec:hct}}

Our testing procedure targets scenarios in which hazard differences may be driven by a small number of time intervals. It follows from previous studies of similar multiple-hypothesis testing situations that a set of individual tests thought to provide the best evidence against the null hypothesis is given by the so-called higher criticism threshold \citep{donoho2008higher,donoho2009feature}, defined as the index $t^\star$ of the P-value maximizing the higher criticism functional $\HC_{t;T}$. Consequently, we define the set of time intervals
\begin{align}
    \label{eq:Delta_def}
\Delta^\star := \{t \, : p_{t} \leq p_{(t^*)}\},\qquad  t^\star = \arg \max_{t \leq T} \HC_{t;T}. 
\end{align}
Figure~\ref{fig:data_intro} illustrates membership in $\Delta^\star$ using an example survival data set. 
    In our situation of comparing survival curves, $\Delta^*$ contains the smallest $t^\star \geq 1$ hypergeometric p-values; these p-values are thought to drive the global difference between the survival curves. The practical value of this identification depends on the context. For instance, time intervals of truly excessive risk suggest points for potential intervention or for conducting further analysis. As such intervals are generally difficult to identify reliably when the departures are small \citep{jin2016rare}, focusing attention on members of $\Delta^*$ likely to increase the utility of such intervention in analogy with feature selection for classification as studied in \cite{donoho2008higher}. To summarize, we propose a single procedure for survival analysis based on $\HCHG$ that combines global testing with the identification of time intervals suspected of excessive risk. 

The multiple-testing approach to survival analysis we promote in this manuscript suggests that feature selection procedures other than the higher criticism threshold of 
\eqref{eq:Delta_def} may also be useful, such as false discovery rate (FDR) controlling  \citep{benjamini1995controlling}. Nevertheless, our theoretical analysis in Section~\ref{sec:analysis} shows that global testing based on FDR controlling is asymptotically powerless in some regimes when our method is still asymptotically powerful. Consequently, in these regimes, HCHG would indicate that the two survival functions of the groups are significantly different at power tending to $1$ as the Type I error tends to $0$, whereas FDR controlling at any false discovery rate parameter would yield a test asymptotically of trivial power. We refer to 
\cite{donoho2008higher,donoho2009feature} for additional discussion about the difference between the two methods in a more general context. 

\subsection{Asymptotic power analysis}
We analyze a test based on $\HCHG$ of \eqref{eq:HC_def} and compare it to other tests using a model involving survival data experiencing non-proportional rare hazard departures. The main purpose of this analysis is to provide a theoretical validation for the method's success compared to existing ones in a framework that emphasizes the rare hazard departure effect we aim to detect. Analysis of this kind is common in mathematical statistics and often leads to useful data analysis tools even under violations of the model's assumptions \citep{donoho2004higher,mukherjee2015hypothesis,arias2015sparse,jin2016rare,pilliat2023optimal}.

Under our framework, the number of at-risk subjects in Group $x$ at time $t$, denoted as $N_x(t)$, experineces a random decay in which the reduction in at-risk subjects at time $t$ follows a Poisson distribution with rate $N_x(t) \bar{\lambda}_t$. Here $\{\bar{\lambda}_t\}_{t=1,\ldots,T}$ is some global baseline hazard sequence hence $\bar{\lambda}_t$ indicates the average time between failure events at the $t$-th interval. The number of at-risk subjects in Group $y$ largely follows the same behavior, except in a few instances in which the rate of the Poisson distribution is $N_y(t) \bar{\lambda}_t'$, where $\bar{\lambda}_t'$ is obtained by perturbing $\bar{\lambda}_t$ upwards. For simplicity, we assumed no censorship. However, it is straightforward to modify the framework to simple right-censorship models in which censoring distribution is independent of mortality events within each group.

We calibrate the number of excessive hazard intervals, the intensity of their departures, and the initial number of subjects to $T$, such that individual perturbations appear on the moderate deviation scale as $T \to \infty$ \citep{RubinSethuraman1965,zeitouni1998large,kipnis2021logchisquared}. Such calibration gives rise to a parameter $r>0$ controlling the intensity of hazard departures and a parameter $\beta \in (0,1)$ controlling their rarity. In this situation, the asymptotic power of $\HCHG$ and other testing procedures experience a phase transition in the following sense. There exists a curve $\left\{\left(\beta, \rho(\beta)\right)\right\}_{\beta \in (0,1)}$ that divides the $(r,\beta)$-plane into two regions. For values $r > \rho(\beta)$, the statistic $\HCHG$ is asymptotically powerful in the sense that there exists a sequence of thresholds for a test based on $\HCHG$ under which the sum of Type I and Type II errors goes to zero. For values of $r < \rho(\beta)$, the sum of Type I and Type II errors goes to one under any sequence of thresholds for a test based on $\HCHG$. The phase transition curve defined by $\rho(\beta)$ turns out to be equal to the phase transition curve of HC in the two-sample normal and Poisson means models under rare and weak perturbations described in \cite{DonohoKipnis2020}. We also discuss the asymptotic properties of additional testing procedures like a test based on Bonferroni's correction (the $\minP$ test) and the false discovery rate controlling procedure. For comparison, our analysis implies that the log-rank test is asymptotically powerless in the entire range of severe rarity $\beta \in (1/2,1)$, whereas HCHG is asymptotically powerful in this range whenever $r > \rho(\beta)$; this situation is illustrated in Figure~\ref{fig:theoretical_PT}. In Section~\ref{sec:simulations}, we exemplify our theoretic derivations numerically by illustrating the Monte-Carlo simulated power of these tests over a grid of configurations of $\beta$ and $r$. The usefulness of our method is not limited to survival data obeying this model, as we demonstrate in Section~\ref{sec:sim_gene_expression} using actual survival data associated with gene expression. 

\begin{figure}
    \centering
    \input{Figs/tkz/PT_th.tex}
    \caption{Phase transition and regions of asymptotic power under the piece-wise homogeneous exponential decay model with rare and weak hazard departures of \eqref{eq:data_model0}-\eqref{eq:lambda_prime_def}. Here $\beta$ controls the number of intervals of excessive hazard and $r$ controls their strength. Our HCHG procedure is asymptotically powerful for $r > \rho(\beta)$. The log-rank test is asymptotically powerless for any $\beta>1/2$. All tests based on randomized hypergeometric tests are asymptotically powerless for $r < \rho(\beta)$.} 
    \label{fig:theoretical_PT}
\end{figure}

While the asymptotic power of $\HCHG$ and other test statistics mirrors previously studied rare and weak signal detection settings \citep{donoho2004higher,cai2014optimal,arias2015sparse,jin2016rare,kipnis2021logchisquared}, our piece-wise exponential decay model introduces additional complexity due to the apparent temporal dependence of events within each group. However, this dependence disappears asymptotically: conditional on the number of at-risk subjects at the start of each interval, the events become independent, and the number of at-risk subjects concentrates around a deterministic value that depends only on the interval and the baseline hazard. A key aspect of our analysis is the calibration of the model's parameters to ensure this concentration, while also guaranteeing that the hypergeometric p-values corresponding to non-null hazard departures exhibit moderately large effects uniformly across time intervals. 
Under these conditions, we establish that the structure of our problem admits an application of the rare and moderate effect detection framework from \cite{DonohoKipnis2020b} and \cite{kipnis2021logchisquared}, from which the main results are derived.

%% file: table_intro_discovery.tex
\begin{tabular}{|c|c|c|c|c|c|c|c|c|c|c|}
\hline
\scriptsize test name: & 
\scriptsize HCHG & 
\scriptsize \begin{tabular}{@{}c@{}}KONP\\(Log-rank)\end{tabular} & 
\scriptsize \begin{tabular}{@{}c@{}}FH\\(0,1)\end{tabular} & 
\scriptsize \begin{tabular}{@{}c@{}}KONP\\(Cauchy)\end{tabular} & 
\scriptsize \begin{tabular}{@{}c@{}}FH\\(1,1)\end{tabular} & 
\scriptsize \begin{tabular}{@{}c@{}}FH\\(0.5,0.5)\end{tabular} & 
\scriptsize Log-rank & 
\scriptsize \begin{tabular}{@{}c@{}}Tarone-\\Ware\end{tabular} & 
\scriptsize \begin{tabular}{@{}c@{}}Gehan-\\Wilcoxon\end{tabular} & 
\scriptsize \begin{tabular}{@{}c@{}}Peto-\\Prentice\end{tabular} \\
\hline
\scriptsize \begin{tabular}{@{}c@{}}proportion of\\true discoveries:\end{tabular} &
\cellcolor{gray!20} 0.66 & 
0.30 & 
0.28 & 
0.27 & 
0.27 & 
0.27 & 
0.27 & 
0.25 & 
0.20 & 
0.20 \\
\hline
\end{tabular}

%% file: table_data_intro.tex
\resizebox{0.4\textwidth}{!}{%
\begin{tabular}{rrrrrr}
\toprule
$t-1$ & $n_x(t-1)$ & $n_y(t-1)$ & $o_x(t)$ & $o_y(t)$ & $p_t$ \\
\midrule
0 & 1000 & 1000 & 11 & 5 & 0.9622 \\
1 & 989 & 995 & 7 & 9 & 0.4062 \\
2 & 982 & 986 & 6 & 6 & 0.6159 \\
3 & 976 & 980 & 2 & 7 & 0.0903 \\
\rowcolor{lightgray} 4 & 974 & 973 & 8 & 20 & 0.0171 \\
5 & 966 & 953 & 6 & 8 & 0.3849 \\
6 & 960 & 945 & 10 & 11 & 0.4854 \\
\rowcolor{lightgray} 7 & 950 & 934 & 8 & 21 & 0.0102 \\
8 & 942 & 913 & 5 & 8 & 0.2705 \\
9 & 937 & 905 & 5 & 5 & 0.6018 \\
\vdots & \vdots & \vdots & \vdots & \vdots \\
13 & 914 & 886 & 6 & 8 & 0.3722 \\
\rowcolor{lightgray} 14 & 907 & 878 & 1 & 15 & 0.0002 \\
15 & 905 & 863 & 6 & 6 & 0.5806 \\
\vdots & \vdots & \vdots & \vdots & \vdots \\
48 & 712 & 668 & 5 & 3 & 0.8343 \\
\rowcolor{lightgray} 49 & 706 & 665 & 6 & 17 & 0.0115 \\
\bottomrule
\end{tabular}}
\resizebox{0.4\textwidth}{!}{%
\begin{tabular}{rrrrrr}
\toprule
$t-1$ & $n_x(t-1)$ & $n_y(t-1)$ & $o_x(t)$ & $o_y(t)$ & $p_t$ \\
\midrule
50 & 700 & 647 & 2 & 3 & 0.4631 \\
51 & 698 & 644 & 7 & 1 & 0.9947 \\
52 & 691 & 643 & 10 & 4 & 0.9621 \\
\rowcolor{lightgray} 53 & 681 & 639 & 7 & 18 & 0.0139 \\
54 & 674 & 621 & 9 & 6 & 0.8100 \\
55 & 665 & 615 & 6 & 2 & 0.9559 \\
56 & 659 & 613 & 3 & 10 & 0.0341 \\
\vdots & \vdots & \vdots & \vdots & \vdots \\
 & & & & \\
76 & 550 & 528 & 3 & 3 & 0.6373 \\
77 & 547 & 525 & 5 & 3 & 0.8425 \\
\rowcolor{lightgray} 78 & 542 & 522 & 2 & 13 & 0.0028 \\
79 & 540 & 509 & 7 & 4 & 0.8678 \\
80 & 533 & 504 & 5 & 4 & 0.7187 \\
81 & 528 & 500 & 8 & 2 & 0.9870 \\
82 & 520 & 498 & 2 & 6 & 0.1299 \\
83 & 518 & 492 & 6 & 3 & 0.8978 \\
84 & 512 & 489 & 0 & 0 & 1.0000 \\
\bottomrule
\end{tabular}
}

%% file: Figs/tkz/PT_th.tex
\pgfdeclarepatternformonly{north east lines wide}%
   {\pgfqpoint{-1pt}{-1pt}}%
   {\pgfqpoint{20pt}{20pt}}%
   {\pgfqpoint{19pt}{19pt}}%
   {
     \pgfsetlinewidth{0.2pt}
     \pgfpathmoveto{\pgfqpoint{0pt}{0pt}}
     \pgfpathlineto{\pgfqpoint{19.1pt}{19.1pt}}
     \pgfusepath{stroke}
    }

\pgfdeclarepatternformonly{north west lines wide}%
   {\pgfqpoint{-1pt}{-1pt}}%
   {\pgfqpoint{20pt}{20pt}}%
   {\pgfqpoint{19pt}{19pt}}%
   {
     \pgfsetlinewidth{0.2pt}
     \pgfpathmoveto{\pgfqpoint{19.1pt}{0pt}}
     \pgfpathlineto{\pgfqpoint{0pt}{19.1pt}}
     \pgfusepath{stroke}
}

\begin{tikzpicture}[scale = 1]
	\begin{axis}[
    width=9.5cm,
    height=5.5cm,
    legend style={at={(0.15,1)},
      anchor=north west, legend columns=1},
    ylabel={$r$ (hazard departure intensity)},
    xlabel={$\beta$ (hazard departure rarity)},
    xtick={.5, .75, .95, 1},
    xticklabels={.5,.75,.95,1},
    ymin=0,
    xmin=0.45,
    xmax=1,
    ymax=2.1,
    title = {Phase transition of tests' power}
    ]

\addplot[color=LRcolor, style=ultra thick, mark=none, 
pattern=north west lines wide,
pattern color=LRcolor
] 
coordinates{(.5,0) (.5,2.5) (1,2.5) (1,-.5)} \closedcycle;

\addlegendentry{\scriptsize $\beta = 1/2$}
  
 
\addplot[domain=0.5:1, mark=none, color=HCcolor, style=ultra thick, 
pattern=north east lines wide,%
pattern color=HCcolor] 
    {2*(x-1/2)*(x<.75) + 2*(sqrt(1)-sqrt(1-x))^2 * (x>=.75) } \closedcycle;

\addlegendentry{\scriptsize $\rho(\beta)$}


\node (PHall) at (axis cs:.87,.39) {};
\node (PHC) at (axis cs:.55,1.3) {};

\end{axis}
\node[right, color=HCcolor, align=center, fill=white] at (PHall) {
\scriptsize all tests are \\ 
\scriptsize asymptotically \\ \scriptsize powerless};
\node[right, color=LRcolor, align=left, fill=white] at (PHC) {
\scriptsize HCHG is asymptotically powerful; \\
\scriptsize log-rank is asymptotically powerless}; 
\end{tikzpicture}

%% file: analysis.tex
\subsection{Exponential Decay Model}
We now introduce a theoretical framework for analyzing the performance of a test based on $\HCHG$ of \eqref{eq:HC_def} and comparing it to other inference procedures. 

Let $x_0$ and $y_0$ be two deterministic constants describing the initial groups' sizes. For $t=0,\ldots, T$, denote by $n_{\star}(t)$ the number of at-risk subjects in Group $\star \in \{x,y\}$. Suppose that there are no censored events and that the sequences $\{O_x(t),O_y(t),N_x(t), N_y(t)\}$ obey
\begin{align}
    N_x(0) = x_0\qquad \text{and} \qquad N_y(0) = y_0,
    \label{eq:data_model0}
\end{align}
and for $t=1,\ldots,T$,
\begin{align}
\begin{cases}
    O_x(t) \sim \Pois( N_x(t-1) \bar{\lambda}_x(t) ) \\
    N_x(t) = \left[N_x(t-1) - O_x(t) \right]^+
\end{cases}
\quad
\text{and}\quad 
    \begin{cases}
    O_y(t) \sim \Pois( N_y(t-1) \bar{\lambda}_y(t)) \\
    N_y(t) = \left[N_y(t-1) - O_y(t)\right]^+,
\end{cases}
    \label{eq:data_model1}
\end{align}
where $[x]^+ = \max\{x,0\}$. The model \eqref{eq:data_model0}-\eqref{eq:data_model1} describes a piece-wise exponential decay of the number of at-risk subjects in either group over time. We consider testing the null hypothesis of identical hazard in both groups
\begin{align}
    H_0 \, :\, \bar{\lambda}_x(t) = \bar{\lambda}_y(t) = \bar{\lambda}_t, \qquad \forall t\in \{1,\ldots,T\}, 
    \label{eq:hyp}
\end{align}
for some unspecified sequence $\{\bar{\lambda}_t\}_{t=1,...,T}$, 
against a situation in which Group $y$ experiences some instances of excessive hazard:
\begin{align}
\label{eq:hyp_perturbation}
    H_1 \,:\,
    \bar{\lambda}_x(t) =\bar{\lambda}_t \quad \text{and} \quad \bar{\lambda}_y(t)  = \sqrt{\bar{\lambda}_t} +  \begin{cases} \sqrt{\delta_t} & t \in I \\
    0 & t \notin I.
    \end{cases}
\end{align}
Here $\delta_t \geq 0$ controls the excess hazard within each interval in the set of non-null intervals $I \subset \{1,\ldots,T\}$. To reflect that we do not know a priori which instances are perturbed, we assume that the membership in $I$ is also random: Each $t$ is included in $I$ with probability $\epsilon$ independently of the other $t$'s. We can write \eqref{eq:data_model0}-\eqref{eq:hyp_perturbation} under the randomness in $I$ as follows. 
\begin{align}
    N_x(0) = x_0 \quad \text{and} \quad N_y(0) = y_0.
    \label{eq:model_full0}
\end{align}
\begin{align}
    H_0\,&:\, \begin{cases}
    O_x(t) \sim  \Pois(N_x(t-1)\bar{\lambda}_t), \quad N_x(t) = \left[N_x(t-1) - O_x(t)\right]^+ \\
    O_y(t) \sim  \Pois(N_y(t-1)\bar{\lambda}_t), \quad
    N_y(t) = \left[N_y(t-1) - O_y(t)\right]^+ 
    \end{cases} \forall t=1,\ldots,T. 
    \label{eq:model_full1}
    \\
    H_1\,&:\, \begin{cases}
    O_x(t) \sim  \Pois(N_x(t-1)\bar{\lambda}_t),\quad 
    N_x(t) = \left[N_x(t-1) - O_x(t)\right]^+ \\
    O_y(t) \sim  (1-\epsilon)\Pois(N_y(t-1)\bar{\lambda}_t)+\epsilon \Pois(N_y(t-1)\bar{\lambda}'_t), \\
    \qquad \qquad N_y(t) = \left[N_y(t-1) - O_y(t)\right]^+ 
    \end{cases} \forall t=1,\ldots,T. 
    \nonumber
\end{align}
where 
\begin{align}
     \bar{\lambda}'_t := \left(\sqrt{\bar{\lambda}_t} + \sqrt{\delta_t} \right)^2. 
     \label{eq:lambda_prime_def}
\end{align}
Namely, the number of events in each group over the $t$-th interval is a random variable that follows a Poisson distribution with expectation proportional to the group's size at the beginning of that interval, unless the number of events exceeds the remaining group's size. Under the null hypothesis, there exists a global unspecified ``base'' failure rate sequence $\{\bar{\lambda}_t\}_{t=1,\ldots,T}$. This sequence governs both groups, $x$ and $y$. Under the alternative, there are roughly $T \cdot \epsilon$ a priori unspecified instances in which the rate of events in Group $y$ is larger than $\bar{\lambda}_t$ by an amount of $\delta_t$ in a square root perturbation \eqref{eq:lambda_prime_def} which is natural in analyzing count data \citep{simpson1987minimum}.

Additional remarks are in order. First, we assume that both $\bar{\lambda}_t$ and $\delta_t$ are very small compared with $T$ and the initial group sizes $x_0$ and $y_0$. Consequently, with probability tending to one, neither group reaches extinction during the study period. Second, the Poisson specification implies that the time between two failure events in Group $x$ within the interval $(t, t+1]$ follows an exponential distribution with mean $\bar{\lambda}_t N_x(t-1)$, truncated at $N_x(t-1)$. We will assume below that this mean is relatively large, ensuring that the total number of events in any interval is also relatively large. Third, the inter-event times in Group $y$ during $(t, t+1]$ approximately follow an exponential distribution with mean $\bar{\lambda}_t N_y(t-1)$ with probability $1-\epsilon$, and with mean $\bar{\lambda}_t' N_y(t-1)$ with probability $\epsilon$. Under our asymptotic setting below, $\bar{\lambda}'_t$ is very close to $\bar{\lambda}_t$, so that the terminal numbers of at-risk subjects $N_x(T)$ and $N_y(T)$ do not, by themselves, distinguish $H_0$ from $H_1$. Finally, one may also formulate the problem in a minimax sense, in which the set $I$ of excessive hazard intervals corresponds to the worst possible choice of $T\epsilon$ intervals. Previous studies in related contexts \citep{chan2017optimal,stoepker2024anomaly} suggest that such a formulation does not affect the asymptotic power of the inference procedures developed below.

Piece-wise exponential decay models as in \eqref{eq:data_model1} are common in survival analysis \citep{feigl1965estimation,friedman1982piecewise}, \citep[Ch. 7]{rodriguez2007lecture}; the departures model \eqref{eq:hyp_perturbation} appears to be new in this context and is analogous to previously-studied high-dimensional heterogeneous detection models involving rare and weak effects \citep{donoho2004higher,hall2008properties,delaigle2009higher,cai2014optimal,donoho2015special,mukherjee2015hypothesis,kipnis2021logchisquared}. As we explain below, a natural calibration of the model's parameters provides a framework for comparing the asymptotic performance of statistical procedures such as $\HCHG$ in this setting.

We anticipate that a test based on $\HCHG$ would exhibit similar properties under variations of \eqref{eq:model_full0}-\eqref{eq:lambda_prime_def} that lead to more complex scenarios, guided by intuition from prior work on higher criticism. While our current model leads to asymptotically vanishing dependencies, in contrast to the persistent dependence structures studied in \cite{hall2008properties,hall2010innovated}, certain extensions, such as those involving elevated hazard across multiple intervals, may induce asymptotically prevailing conditional dependencies across time. Addressing such cases may benefit from techniques developed in these works, and we leave these extensions as future work.

\subsection{Calibration}
We consider an asymptotic setting in which $T\to \infty$, while the other parameters 
$\epsilon$, $\delta_t$, $x_0$, $y_0$ and $\bar{\lambda}_t$ are calibrated to $T$ and define a sequence of local alternatives $H_1^{(T)}$ to $H_0^{(T)}$. We introduce additional parameters $r$ and $\beta$ to describe our calibration of \eqref{eq:model_full0}-\eqref{eq:lambda_prime_def}, as summarized in Table~\ref{tbl:calibration}. Our calibration is chosen such that non-null hypergeometric p-values in \eqref{eq:pvals_def} correspond to a rare moderate departure setting which arises when the Poisson rates of affected intervals deviate uniformly on the moderate scale \citep{kipnis2021logchisquared}. Calibrating the model in other ways may lead to different asymptotic power behavior of some inference procedures, in analogy with other rare and weak effect models with non-moderate departures 
\citep{arias2015sparse,jin2016rare,arias2017distribution,DonohoKipnis2020,DonohoKipnis2020b}.
\begin{table}[ht]
    \begin{center}
    \begin{tabular}{|c||c|p{7cm}|p{2cm}|}
    \toprule
         parameter & reference & description & calibrating parameter \\
         \midrule
         \midrule
         $\epsilon$ & \eqref{eq:calibration_eps} & proportion of non-null hazard departures  & $\beta \in (0,1)$ 
          \\
         \midrule 
         $\delta_t$ & \eqref{eq:calibration_delta} & departure size (Hellinger shift in the \newline hazard of non-null occurrences)  & $r \geq 0$  \\
         \midrule 
         $\bar{\lambda}_t$ & \eqref{eq:calibration_rates} & baseline hazard sequence &  \\
         \bottomrule
    \end{tabular}
    \end{center}
    \caption{Parameters of the piece-wise exponential decay under rare and weak hazard departures theoretical framework.}
    \label{tbl:calibration}
\end{table}

We assume that the initial group sizes $x_0$ and $y_0$ go to infinity as $T$ goes to infinity,
while they are asymptotically equivalent in the sense that
\begin{align}
\label{eq:calibration_kappa}
\frac{x_0}{y_0} \to 1,\quad \text{as}\quad T \to \infty.
\end{align}
Additionally, their increase is limited such that
\begin{align}
\label{eq:calibration_initial}
 \lim_{T \to \infty} \frac{x_0}{T^{1+a}} = 0, \qquad \forall a>0.
\end{align}
We calibrate the rarity parameter $\epsilon$ to $T$ according to 
\begin{align}
\label{eq:calibration_eps}
    \epsilon := \epsilon(T) = T^{-\beta},\quad \beta \in (1/2,1).
\end{align}
The complementary situation of $\beta < 1/2$ leads to non-rare asymptotic power behavior under moderate departures, as is well understood from other rare and weak effect studies \citep{jin2003detecting,arias2015sparse}. 
We calibrate the effect size parameter $\delta$ to $T$ according to
\begin{align}
\label{eq:calibration_delta}
    \delta_t : = \delta(t,T) = \frac{r}{2}\frac{ \log(T)}{n(t)},\qquad n(t) :=  \frac{x_0+y_0}{2} e^{-\sum_{s \leq t}\bar{\lambda}_s}.
\end{align}
We summarize the descriptions of model parameters $\epsilon$, $\delta_t$ and $\bar{\lambda}_t$ and their connection to $r$ and $\beta$ in Table~\ref{tbl:calibration}. 

As $T\to \infty$, a standard concentration argument implies $\frac{N_x(t)}{n(t)} \sim \frac{N_y(t)}{n(t)} \to 1$ in probability uniformly in $t\leq T$; see Lemma~4 in the supplement \cite{survival2023supp}. Namely, $n(t)$ is approximately the number of at-risk subjects at time $t$ in either group. This explains the effect size calibration \eqref{eq:calibration_delta} as a departure relative to the number of at-risk subjects. We further assume
\begin{align}
\label{eq:calibration_rates}
\min_{t \leq T} \frac{\bar{\lambda}_t x_0} {\log(T)} \to \infty \quad \text{while} \quad \max_{t \leq T} \bar{\lambda}_t T \leq M, 
\end{align}
for some finite $M$ that is independent of $T$. Hence, the decay rates $\bar{\lambda}_t$ asymptotically vanish, but not too rapidly.

Note that \eqref{eq:calibration_rates} implies 
\[
\ex{O_x(t)} \approx \ex{\bar{\lambda}_t N_x(t-1)} = \bar{\lambda}_t x_0(1-\bar{\lambda}_t)^t \approx  
\bar{\lambda}_t x_0 e^{- \sum_{s \leq t} \bar{\lambda}_s } \geq 
\bar{\lambda}_t x_0 e^{-M} \to \infty,
\]
hence $\ex{O_x(t)} \to \infty$ and likewise $\ex{O_y(t)} \to \infty$ at rates faster than $\log(T)$. Furthermore, conditions \eqref{eq:calibration_kappa}-\eqref{eq:calibration_rates} ensure that the expected proportion relative to the initial size of at-risk subjects in each group at interval $t$ converge in probability to $\exp\{-\sum_{s \leq t} \bar{\lambda}_s\}$ (see Lemma~4 in the supplement \cite{survival2023supp}). In particular, $N_y(t)/(N_x(t) + N_y(t))$ remains roughly constant in $t$ at around $1/2$, and it is generally impossible to recognize any difference in the failure rates by considering this ratio at any given $t=1,\ldots,T$.

\subsection{Asymptotic power and phase transition}
A statistic $U_T$ based on the data $\{N_x(t), N_y(t), O_x(t), O_y(t)\}_{t=0}^T$ under the setting \eqref{eq:hyp} is said to be asymptotically \emph{powerful} if there exists a sequence of thresholds $\{h(T)\}_{T=1,2,\ldots,}$ such that 
\begin{align*}
    \Prp{U_T \geq h(T) \mid H_0 } + \Prp{U_T < h(T) \mid H_1 } \to 0.
\end{align*}
Conversely, $U_T$ is said to be asymptotically \emph{powerless} if 
\begin{align*}
    \Prp{U_T \geq h(T) \mid H_0 } + \Prp{U_T < h(T) \mid H_1 } \to 1,
\end{align*}
for any sequence of threshold $\{h(T)\}_{T=1,2,\ldots,}$. In words, the asymptotic powerfulness of $U_T$ means that a sequence of tests for \eqref{eq:hyp} based on $U_T$ approaching full power exists, whereas asymptotic powerlessness says that any sequence of tests based on $U_T$ asymptotically has trivial power. See the references \cite{donoho2004higher,arias2005near,donoho2015special} for additional discussions of these definitions.

\subsection{Asymptotic power of HCHG}

Under the setting \eqref{eq:model_full0}-\eqref{eq:lambda_prime_def} and the calibration \eqref{eq:calibration_kappa}-\eqref{eq:calibration_rates}, Theorem~\ref{thm:HC_powerful} below shows that the curve 
\begin{align}
\label{eq:rho}
\rho(\beta) &:=
\begin{cases}
    2(\beta - 1/2) & \frac{1}{2} < \beta <  \frac{3}{4}, \\
    2\left(1-\sqrt{1 -\beta }\right)^2 & \frac{3}{4} \leq  \beta < 1,
    \end{cases}
\end{align}
characterizes a region in the parameter space $(\beta,r)$ in which $\HCHG$ is asymptotically powerful; see an illustration in Figure~\ref{fig:theoretical_PT}. 
\begin{theorem}
    \label{thm:HC_powerful}
    Consider testing $H_0$ versus $H_1$ as in \eqref{eq:model_full1} when $x_0$, $y_0$, $\{\bar{\lambda}_t\}$, $\epsilon$, and $\delta$ are calibrated to $T$ as in \eqref{eq:calibration_kappa}-\eqref{eq:calibration_rates}. If $r > \rho(\beta)$, 
    $\HCHG$ of \eqref{eq:HC_def} is asymptotically powerful. 
\end{theorem}

A statement about the ineffectiveness of $\HCHG$ and other test statistics is provided in Theorem~\ref{thm:HC_powerless} below. This statement focuses on p-values obtained by randomizing the hypergeometric tests of \eqref{eq:pvals_def} such that each P-value has a continuous distribution yet dominated by its non-randomized version. For example, replace \eqref{eq:pvals_def} by
\begin{align}
\label{eq:pval_randomized}
\tilde{\pi}(x, y;  n_x, n_y) & = \Prp{\HG(n_x + n_y, n_y, x+y) \geq y} \\
& \quad - U \cdot \Prp{\HG(n_x + n_y, n_x, x+y) = y}\nonumber,
\end{align}
where $U$ is uniformly distributed over $(0,1)$ and independent between tests. Such randomization is common in decision theory \cite[p. 101]{cox1979theoretical} and typically improves the power of multiple testing procedures. From an information theoretic perspective, randomization is necessary for meaningful comparison of experiments and statements about the impossibility of inference \citep{le2012asymptotic,brown2002asymptotic,nussbaum2006constructive}; see a related discussion in a similar discrete two-sample setting in
\cite{DonohoKipnis2020}.
\begin{theorem}
    \label{thm:HC_powerless}
    Consider testing $H_0$ versus $H_1$ as in \eqref{eq:model_full1} when $x_0$, $y_0$, $\{\bar{\lambda}_t\}$, $\epsilon$, and $\delta$ are calibrated to $T$ as in \eqref{eq:calibration_kappa}-\eqref{eq:calibration_rates}. Let $\tilde{p}_1,\ldots,\tilde{p}_T$ be p-values obtained from the randomized hypergeometric tests \eqref{eq:pval_randomized}. If $r < \rho(\beta)$, all tests based on $\tilde{p}_1,\ldots,\tilde{p}_T$ are asymptotically powerless. 
\end{theorem}

The proofs of Theorems~\ref{thm:HC_powerful} and \ref{thm:HC_powerless} (in the supplement \citep{survival2023supp}) rely on known properties of the asymptotic power of higher criticism of 
rare moderately departed p-values from \cite{kipnis2021logchisquared}. The proof builds on a series of technical results showing that the hypergeometric p-values of \eqref{eq:pvals_def} under the hypothesis testing setting \eqref{eq:model_full0}-\eqref{eq:lambda_prime_def} and the calibration \eqref{eq:calibration_kappa}-\eqref{eq:calibration_rates} correspond to the rare moderate departures (RMD) model in the following sense. Define $\alpha(q,\rho):= (\sqrt{q}-\sqrt{\rho})^2$. 
In the supplement \citep{survival2023supp}, we show that the p-values of \eqref{eq:pvals_def} obey
\begin{align}
\label{eq:pvalue_tail}
    \lim_{T \to \infty} \max_{t=1,\ldots,T} \left| \frac{-2\log(\Prp{ -2\log(p_t) \geq 2q\log(T)\mid t \in I })}{\log(T)} - \alpha(q,r/2) \right| = 0,
\end{align}
for $q > r/2 > 0$. The limit in \eqref{eq:pvalue_tail} says that at any interval $t$ of truly elevated hazard (indicated by $t\in I$), the tail of $-2\log(p_t)$ on the moderate deviations scale behaves as the tail of a non-central chisquared random variable over one degree of freedom. Namely, the distribution of $p_t$ conditioned on $t \in I$ satisfies
\[
-2\log(p_t) \overset{D}{\approx} \left(\sqrt{r \log(T)} + Z \right)^2,\qquad Z\sim \Ncal(0,1),
\]
with the approximation in the sense of \eqref{eq:pvalue_tail}. Furthermore, $p_t$ is independent of previous $p_s$ for $s<t$ conditioned on the number of at-risk subjects in each group $N_x(t)$ and $N_y(t)$. Since these random variables concentrate around $x_0\exp\{-\sum_{s \leq t} \bar{\lambda}_s\}$ and $y_0\exp\{-\sum_{s \leq t} \bar{\lambda}_s\}$, respectively (see Lemma~4 in the supplement \cite{survival2023supp}), the asymptotic joint distribution of the p-values converges to a product distribution so the setting of \cite{kipnis2021logchisquared} applies.

\subsection{Asymptotic powerlessness of the log-rank test}
The Log-Rank test is defined as follows. 
Set $n(t):= n_x(t) + n_y(t)$, $o(t):=o_x(t) + o_y(t)$, and  
\begin{align*}
    e_t & := \frac{n_y(t-1)}{n(t-1)} \left(o_x(t)+o_y(t)\right),\\
    v_t & :=\frac{n_y(t-1)n_x(t-1)}{n(t-1)-1} \frac{\left(o_x(t)+o_y(t)\right)}{n(t-1)} \left(1 - \frac{\left(o_x(t)+o_y(t)\right)}{n(t-1)} \right). 
\end{align*}
The log-rank test statistic is
\begin{align}
    \label{eq:logrank}
\LR := \frac{\sum_{t=1}^T o_y(t) - \sum_{t=1}^T e_t}{\sqrt{\sum_{t=1}^T v_t}}
\end{align}
and we reject for large values of $\LR$ \citep{mantel1966evaluation,cox1975partial}. 

\begin{theorem}
    \label{thm:LR_powerless}
    Consider testing $H_0$ versus $H_1$ as in \eqref{eq:model_full1} when $x_0$, $y_0$, $\bar{\lambda}$, $\epsilon$, and $\delta$ are calibrated to $T$ as in \eqref{eq:calibration_kappa}-\eqref{eq:calibration_rates}. $\LR$ is asymptotically powerless. 
\end{theorem}
    The proof of Theorem~\ref{thm:LR_powerless} (in the supplement \cite{survival2023supp}) is based on the asymptotic normality of $\LR$ and the analysis of its first two moments under either hypothesis.
We note that the asymptotic behavior of the log-rank statistic appears to be analogous to that of the Fisher combination statistics 
\[
F_T := 2\sum_{t=1}^T \log(1/p_t)
\]
in general rare moderate departures models  \citep{kipnis2021logchisquared}. Figure~\ref{fig:theoretical_PT} illustrates the region in which $\LR$ is asymptotically powerless. 

\subsection{Asymptotic power of other multiple testing procedures}

Asymptotic characterizations of several multiple-testing procedures involving p-values obeying the rare moderate departures formulations are available in \cite{kipnis2021logchisquared}. These include the $\minP$ that is based on $p_{(1)}$ and a test that is based on Benjamini-Hochberg's false discovery rate (FDR) functional
\begin{align*}
\mathrm{FDR}^*(p_1,\ldots,p_T) :=  \min_{1\leq t \leq T} \frac{p_{(t)}}{t}. 
\end{align*}
Both tests have the same phase transition curve separating the region of asymptotic powerfulness from powerlessness, given by
\begin{align*}
    \rho_{\minP}(\beta) :=
    2\left(1 - \sqrt{1 -\beta }\right)^2,\qquad  1/2 < \beta < 1.
\end{align*}
Namely, whenever $3/4 < \beta < 1$, a test based either on the minimal P-value or on $\mathrm{FDR}^*$ are asymptotically powerful in the same region in which $\HCHG$ is asymptotically powerful. On the other hand, when $1/2< \beta < 3/4$, there exists a region in which $\HCHG$ is asymptotically powerful but these other two tests are not. This situation is analogous to other two-sample multiple testing settings under rare and moderately large departures \citep{DonohoKipnis2020}.

We summarize the regions in which different tests are asymptotically powerful or powerless in Table~\ref{tbl:power}. 

\begin{table}[ht]
    \begin{center}
    \begin{tabular}{|l|c|l|}
    \toprule
        method & test statistic & asymptotic power\\
        \midrule
        \midrule
        Higher criticism & $\HCHG$ & 
        powerful when $\beta \in (1/2,1)$,  $r > \rho(\beta)$ \\
        \midrule
        Bonferroni ($\minP$) & 
       $1 / p_{(1)} $  & powerful when $\beta \in (1/2,1)$, $r > \rho_{\minP}(\beta)$ \\
       \midrule
       False discovery rate &  $\max_{t} \frac{t}{p_{(t)}} $ & powerful when $\beta \in (1/2,1)$, $r > \rho_{\minP}(\beta)$ \\
       \midrule
       Log-rank & $\LR$ & powerless when $\beta > 1/2$\\
       \bottomrule
    \end{tabular}
    \end{center}
    \caption{Region of asymptotic power of various testing methods under the two-sample Poisson decay model with rare and weak hazard departures \eqref{eq:model_full1}. 
Higher criticism of hypergeometric p-values ($\HCHG$) has the largest region of powerfulness among the tests in the table. }
    \label{tbl:power}
\end{table}

%% file: simulations.tex
\newcommand{\testSizeDense}{0.55}
\newcommand{\testSizeSparse}{0.6}
\newcommand{\sigLevel}{0.05}
\newcommand{\npow}{5}
\newcommand{\xiSparse}{0.8}
\newcommand{\xiDense}{1.4}
\newcommand{\nMonte}{1,000}
\newcommand{\nMonteNull}{{100,000}}
\newcommand{\nMonteII}{{50,000}}
\newcommand{\numberPatients}{{3,069}}
\newcommand{\numberGenes}{{8,702}}
\newcommand{\numberGenesRemovedI}{{557}}
\newcommand{\numberGenesRemovedII}{{4674}}
\newcommand{\numberTested}{{4028}}
\newcommand{\numberGenesFull}{{9,259}}
\newcommand{\numberIntervals}{{82}}

\newcommand{\numDiscoveriesOOs}{1671}
\newcommand{\numDiscoveriesOIs}{1451}
\newcommand{\numDiscoveriesIOs}{163}
\newcommand{\numDiscoveriesIIs}{743}
\newcommand{\numDiscoveriesOO}{1669}
\newcommand{\numDiscoveriesOI}{1450}
\newcommand{\numDiscoveriesIO}{165}
\newcommand{\numDiscoveriesII}{744}
        

\subsection{Simulated Data and Empirical Phase Transition}

We conduct a sequence of Monte-Carlo experiments with $x_0 = y_0 = T \log(T)$ and the baseline rate $\bar{\lambda} = 2/T$. Similarly to \cite{DonohoKipnis2020}, we consider points $(\beta,r)$ in a grid $I_r \times I_\beta$ covering the range $I_r \subset [0,2.1]$, $I_\beta \subset (0.45,0.95)$. For each test statistic $U$, we first find the $1-\alpha$ empirical quantile of $U$ under the null hypothesis $r=0$ using $N_0 = \nMonteNull$ Monte-Carlo experiments. We denote this quantile as $\hat{q}^{1-\alpha}(U)$. Next, for each configuration $(\beta,r)$, we conduct $N=\nMonte$ Monte-Carlo experiments and consider the number of instances in which $U$ exceeds $\hat{q}^{1-\alpha}(U)$. We denote this number by $\hat{B}(U,\alpha,\beta,r)$. We declare that $\hat{B}(U,\alpha,\beta,r)$ is \emph{substantial} if we can reject the hypothesis $H_{0,\alpha}~~:~~ \hat{B}(U,\alpha,\beta,r) \sim \Bin(N,\alpha)$ at level $\alpha_1$. Namely, $\hat{B}(T,\alpha,\beta,r)$ is substantial if $\Pr \left( \Bin(N,\alpha) \geq \hat{B}(T,\alpha,\beta,r) \right) \leq \alpha_1$. Next, we fix $\beta \in I_\beta$ and focus on the strip $\{(\beta,r)\}_{r \in I_r}$. We construct the binary-valued vector indicating those values of $r$ for which $\hat{B}(T,\alpha,\beta,r)$ is substantial. To this vector, we fit a logistic response model. The phase transition point of the strip $\{(\beta,r),\, r \in I_r\}$ is defined as the point $r=\hat{\rho}(\beta)$ at which the fitted response equals $0.5$. The empirical phase transition curve is defined as $\{ \hat{\rho}(\beta)\}_{\beta \in I_{\beta}}$. 

The top panels of Figure~\ref{fig:sim} illustrate the Monte-Carlo simulated power $\hat{B}(T,\alpha,\beta,r)/N$ and the empirical phase transition curves of $\HCHG$ and $\LR$ along with their theoretical counterparts. The results illustrated in these figures support our theoretical finding in Theorems~\ref{thm:HC_powerful} and \ref{thm:HC_powerless}, establishing $\rho(\beta)$ of \eqref{eq:rho} as the boundary between the region where HCHG has asymptotically maximal power and the region where it has asymptotically no power. We also show the Monte-Carlo simulated power and empirical phase transition of the log-rank test; the region of powerfulness is smaller than that of $\HCHG$, in agreement with our theoretical result in Theorem~\ref{thm:LR_powerless}. The empirical phase transitions of some weighted versions of the log-rank test we experimented with are similar but typically inferior to that of the log-rank. These statistics include: Tarone-Ware \citep{tarone1977distribution}, Gehan-Wilcoxon
\citep{gehan1965generalized}, and Fleming-Harrington with $(p,q) \in \{(1,0),(0.5,0.5),(1,1)\}$ \citep{harrington1982class}. The bottom panel in Figure~\ref{fig:sim} indicates configurations on the grid $I_\beta \times I_r$ in which the empirical power of a test based on $\HCHG$ at the level $\sigLevel$ is significantly better or worse than a test based on different statistics. 

\begin{figure}[ht!]
\begin{center}
\input{Figs/tkz/PT_sim_HC.tex}
\caption{Empirical power of higher criticism of hyper-geometric (HCHG) p-values (left) and log-rank (right) at level $\alpha=0.05$. The curves $\rho(\beta)$ and $\hat{\rho}(\beta)$ are the theoretical and Monte-Carlo simulated phase transitions of HCHG, respectively. The line $\beta=0.5$ and the curve $\hat{\rho}_{\LR}(\beta)$ are the theoretical and Monte-Carlo simulated phase transition of the log-rank statistics of \eqref{eq:logrank}, respectively. 
}
\label{fig:sim}
    \end{center}
\end{figure}

\begin{figure}
    \begin{center}
        \input{power_diff.tex}
    \end{center}
    \caption{Configurations with significant empirical power differences between a test based on 
    $\HCHG$ of \eqref{eq:HC_def} and tests based on other statistics. A gray point indicates no significant power difference towards any statistic. We used $N=1,000$ experiments in each $(r,\beta)$ configuration. Each experiment simulates a sample from the piece-wise exponential decay model \eqref{eq:data_model1} with rare and weak departures over $T=1,000$ time intervals. }
\end{figure}

\subsection{Demonstration for gene expression data}
\label{sec:sim_gene_expression}

\subsubsection{Setup}
The SCANB dataset \citep{saal2015sweden} records mortality events over time of $\numberPatients$ breast cancer patients. It also includes the expression level of
$\numberGenesFull$ genes in each patient. We removed $\numberGenesRemovedI$ genes whose response contains repeated values. For each gene $g$ of the remaining $\numberGenes$, we divide the patients into two groups: Group $x$ consists of patients whose expression for $g$ is at or below the median value for $g$, and Group $y$ contains all patients whose expression is larger than this median. This process partitions the patients into two groups of roughly equal sizes, denoted by $x_0(g)$ and $y_0(g)$, $\numberGenes$ partitions overall. In each partition, we consolidated all events into $T=\numberIntervals$ intervals of approximately 28 days each. 

\subsubsection{Simulating null distribution}
Over $N=\nMonteII$ iterations, we randomly assign half of the patients to Group $x$ and the other half to Group $y$. This assignment leaves the original correspondence between censorship and event times but removes group associations that, in the actual data, are driven by biology (gene expression). The empirical quantile resulting from this permutation procedure is known to be useful for level testing when the censoring distributions within groups are identical \cite[Sec. 5]{heimann1998permutational}. Consequently, we tested all genes for equality of the censoring distributions using a Kolmogorov Smirnoff test; we removed $\numberGenesRemovedII$ genes in which the difference is significant at level $0.05$. 
We evaluated the $0.95$ empirical quantile of the size-$N$ sample of values of $\HCHG$ of \eqref{eq:HC_def}, denoted $\hat{q}_{0}^{0.95}(\HCHG)$ using the randomization procedure described above. The full histogram of $\HCHG$ values is provided in Figure~\ref{fig:histogram}. 
The difference between the simulated null values of $\HCHG$ under random group assignments and uniformly sampled p-values, e.g., as reported in \cite{donoho2004higher}, follows from the super uniformity of the hypergeometric p-values. Consequently, the Z-scores in the HC calculation \eqref{eq:HC_def} are biased downwards hence their maximum is also much smaller and may even be negative. 

\begin{figure}[ht!]
    \centering
    \input{Figs/tkz/histogram}
    \caption{Simulated null distribution of $\HCHG$. Histogram of $\HCHG$ over $N=\nMonteII$ random group assignments with event times taken from the SCANB gene expression dataset. The $0.95$-th quantile is indicated. 
    }
    \label{fig:histogram}
\end{figure}

\subsubsection{Testing} 
For each gene $g$ of the remaining $\numberTested$, we applied our testing procedure in Algorithms~\ref{alg:1} and \ref{alg:2} to check for an increased hazard in either group. Namely, we report the existence of any effect associated with $g$ if $\HCHG(g)$ exceeds $\hat{q}^{0.95}_0(\HCHG)$ or if $\mathrm{Rev}[\HCHG(g)]$ exceeds $\hat{q}^{0.95}_0(\HCHG)$, where 
 $\mathrm{Rev}[\HCHG(g)]$ is obtained from Algorithm~\ref{alg:1} after switching the roles of the $x$-series and the $y$-series. We report the existence of a strictly one-sided effect associated with $g$ if $\HCHG(g)$ exceeds $\hat{q}^{0.95}_0(\HCHG)$ while 
 $\mathrm{Rev}[\HCHG(g)]$ does not exceed $\hat{q}^{0.95}_0(\HCHG)$, or vice versa. We also used the log-rank test based on 
 $\hat{q}^{0.95}_0(\LR)$, the simulated $0.95$ quantile of the log-rank statistic $\LR$ of \eqref{eq:logrank} under $H_0$, as well as several weighted versions of the log-rank test proposed in the literature to discover non-proportional hazard departures \citep{yang2010improved}.

\input{Figs/tkz/venn_1sided.tex}

\subsubsection{Results}
We report the number of genes found significant at the level $\sigLevel$ by each testing method in Figure~\ref{fig:venn_2sided}. For example, in testing for a strictly one-sided effect, HCHG identified $\numDiscoveriesIOs$ significant genes that the log-rank test did not report as significant. We list some of these genes in Table~\ref{tbl:results_SCANB}, where we also report on the empirical p-values with respect to each test statistic associated with these genes. The Kaplan-Meier survival curves corresponding to three example genes are illustrated in Figure~\ref{fig:real_data}. In these figures, we also indicate time intervals driving change between the groups according to the higher criticism thresholding procedure of \eqref{eq:Delta_def}.
\input{table_results_SCANB}

\input{Figs/tkz/fig_real_data.tex}

%% file: Figs/tkz/PT_sim_HC.tex
\newcommand{\level}{.05}

\begin{tikzpicture}[scale = 1]
	\begin{axis}[
    width=6.3cm,
    height=5cm,
    legend style={at={(1,1)},
      anchor=north east, legend columns=1},
    ylabel={$r$ (intensity)},
    xlabel={$\beta$ (rarity)},
    xtick={0.5,0.75,.95},
    ymin=0,
    xmin=0.45,
    xmax=1,
    ymax=2.1,
    ]

\addplot[domain=.4:.43, color=HCcolor, style=dashed, line width= 1pt] {x};
\addlegendentry{\scriptsize $\hat{\rho}(\beta)$}

\addplot[domain=0.75:1, color=HCcolor, style= ultra thick] 
    {2*(sqrt(1)-sqrt(1-x))^2};
    \addlegendentry{\scriptsize $\rho(\beta)$}

\addplot graphics[xmin=0.425,xmax=0.975,ymin=0,ymax=2.3,
		includegraphics={trim=75 57 145 40, clip, scale=.25}]{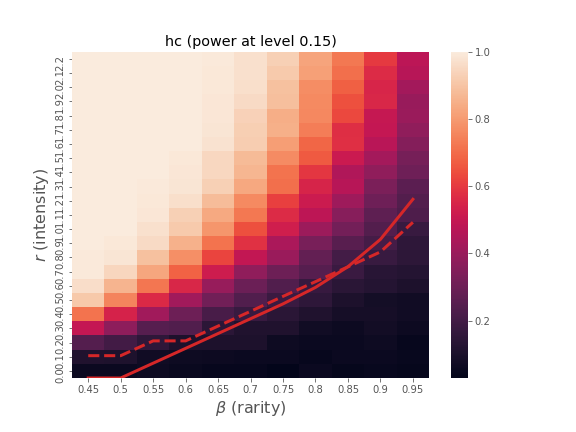};
	
\addplot[domain=0.5:0.75, color=HCcolor, style=ultra thick] 
    {2*(x-1/2)};

\addplot[domain=0.75:1, color=HCcolor, style=ultra thick] 
    {2*(sqrt(1)-sqrt(1-x))^2};
    
\end{axis}
\end{tikzpicture}
\begin{tikzpicture}[scale = 1]
	\begin{axis}[
    width=6.3cm,
    height=5cm,
    legend style={at={(1,1)},
      anchor=north east, legend columns=1},
    xlabel={$\beta$ (rarity)},
    xtick={0.5,0.75,.95},
    yticklabels={},
    ymin=0,
    xmin=0.45,
    xmax=1,
    ymax=2.1,
    ]

\addplot[domain=.4:.43, color=LRcolor, line width= 1pt] {x};
\addlegendentry{\scriptsize $\beta=0.5$}

\addplot[domain=.4:.43, color=LRcolor, style=dashed, line width= 1pt] {x};
\addlegendentry{\scriptsize $\hat{\rho}_{\LR}(\beta)$}

\addplot[domain=.4:.43, color=HCcolor, style=ultra thick] {x};
\addlegendentry{\scriptsize $\rho(\beta)$}


\addplot graphics[xmin=0.425,xmax=0.975,ymin=0,ymax=2.3,
		includegraphics={trim=75 57 145 40, clip, scale=.25}]{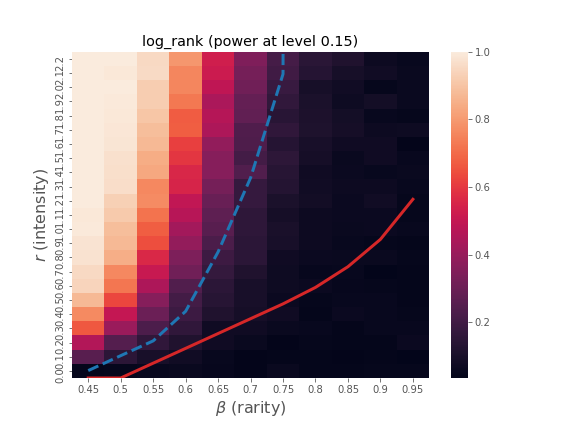};

\addplot[color=LRcolor, style= ultra thick] 
    coordinates{(.5, 0) (0.5, 2.3)};
 
\addplot[domain=0.5:0.75, color=HCcolor, style=ultra thick, opacity=1] 
     {2*(x-1/2)};

\addplot[domain=0.75:1, color=HCcolor, style=ultra thick, opacity=1] 
     {2*(sqrt(1)-sqrt(1-x))^2};
    
\end{axis}
\node at (-0.5, 1.75) 
{\includegraphics[scale=.3] {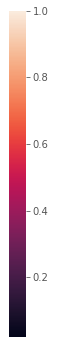}};
\end{tikzpicture}

%% file: power_diff.tex
\includegraphics[scale=.4]{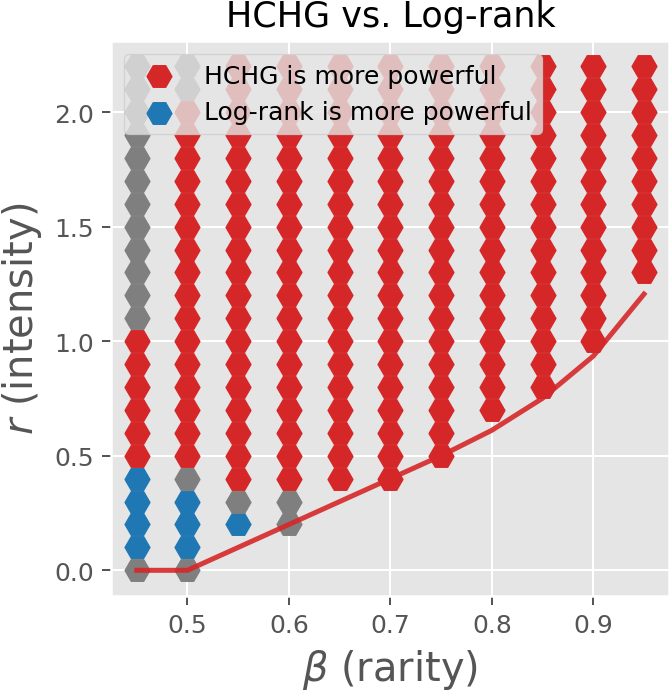}
\includegraphics[scale=.4]{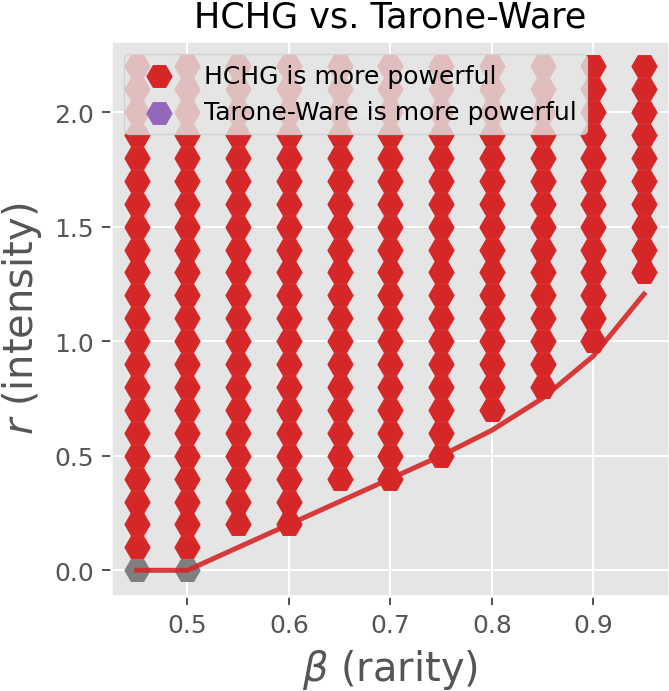}
\includegraphics[scale=.4]{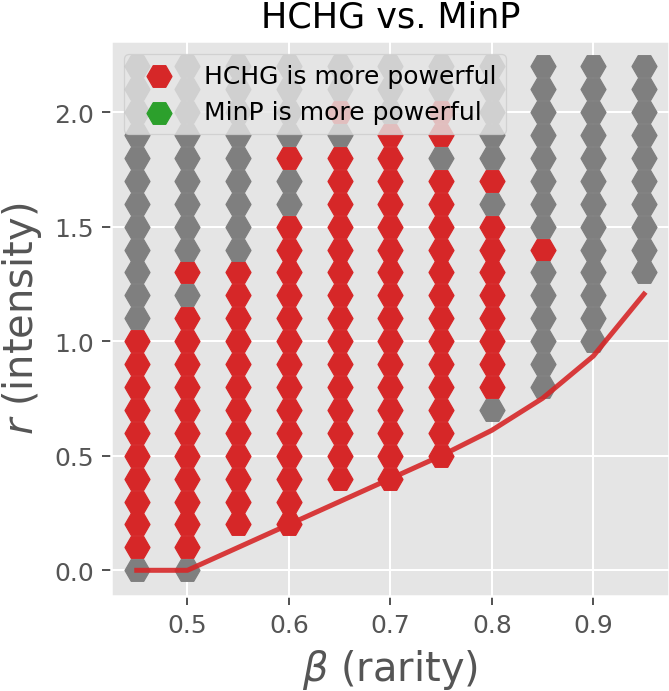}

%% file: Figs/tkz/histogram.tex
\begin{tikzpicture}[scale = 1]
	\begin{axis}[
    width=7.5cm,
    height=3cm,
    legend style={at={(.91,1)},
      anchor=north east, legend columns=1},
    ylabel={\scriptsize frequency},
    xlabel={\scriptsize $\HCHG$ values},
    xtick={-5,-4,-3,-2,-1,0,1,2},
    xticklabels={,,,\scriptsize -2,, \scriptsize 0 , ,\scriptsize 2},
    ytick={0,0.8},
    yticklabels={,\scriptsize 0.8},
    ymin=0,
    xmin=-5,
    xmax=2,
    ymax=0.8,
    ]
\addplot graphics[xmin=-6,xmax=2,ymin=0,ymax=0.8,
		includegraphics={ trim=80 17 0 0, clip, scale=.5}]{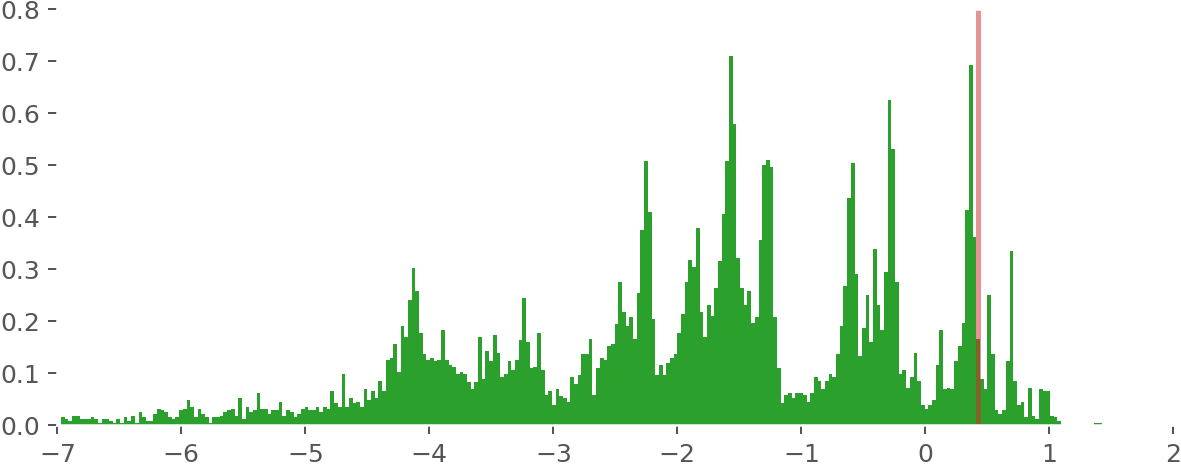};


\node (critval05) at (axis cs:0.65,.67) {};

    
\end{axis}
%



\end{tikzpicture}

%% file: Figs/tkz/venn_1sided.tex
\begin{figure}[!ht]
    \centering
    \includegraphics[scale=.29]
    {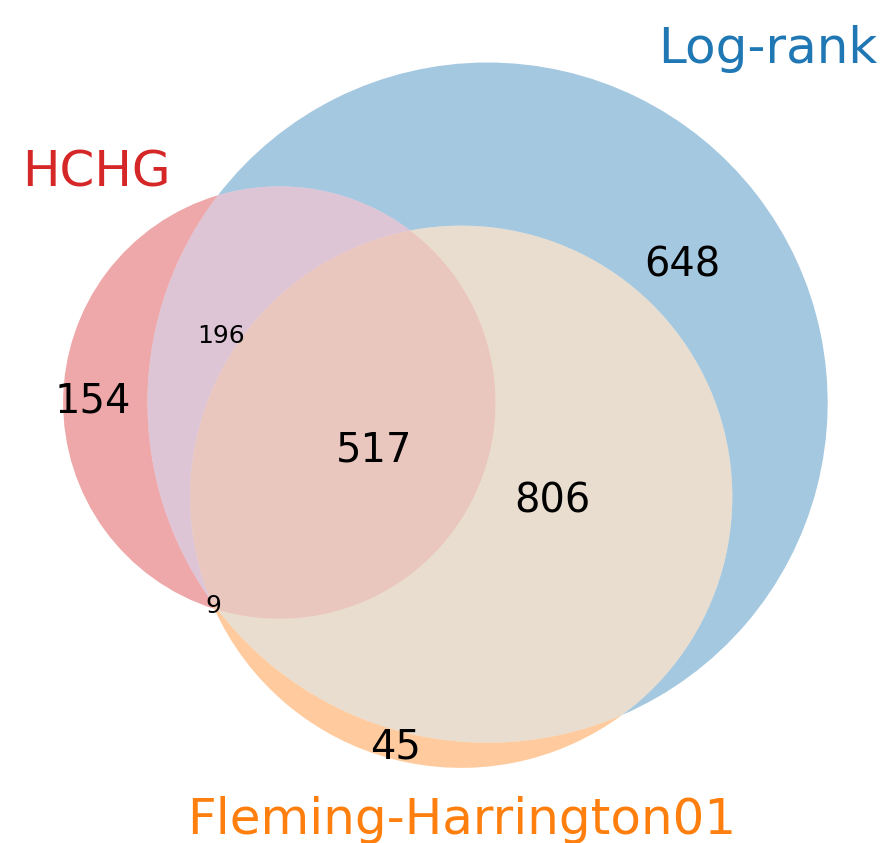}
    \includegraphics[scale=.29]
    {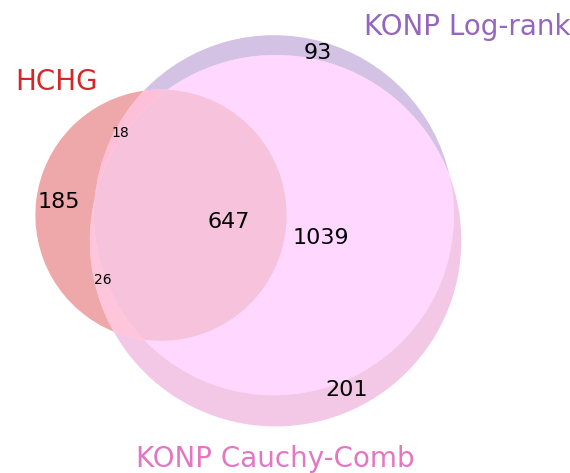}
    \includegraphics[scale=.29]
    {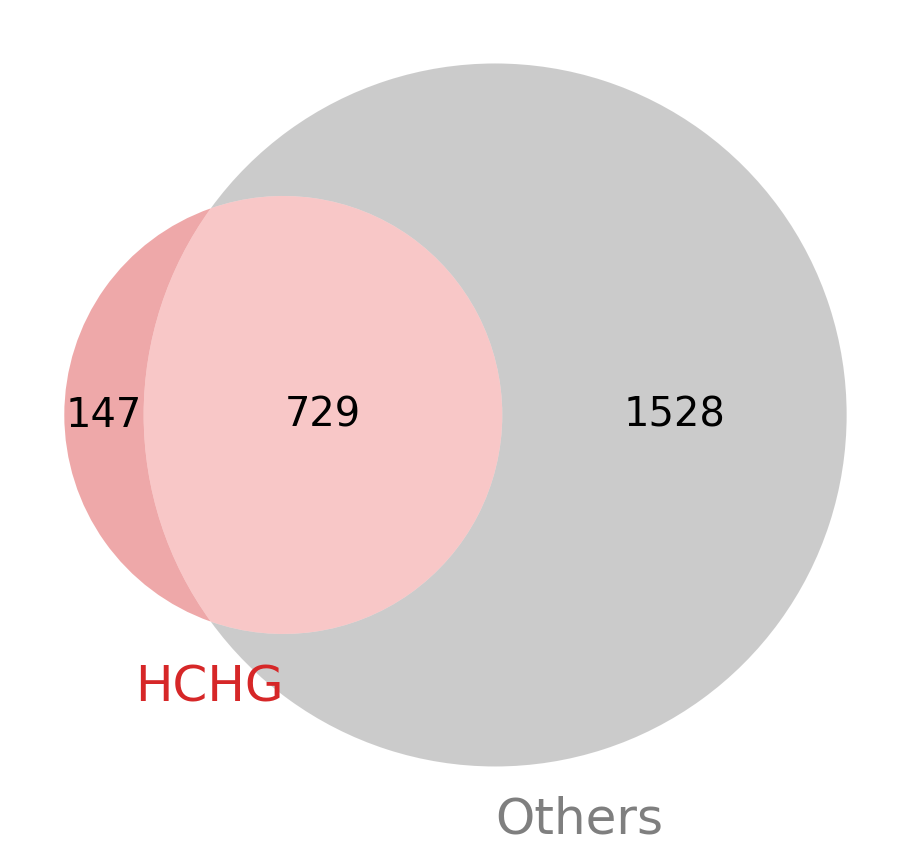}
    \caption{Number of genes with expression levels significantly ($\alpha=\sigLevel$) associated with survival according to the higher criticism of hypergeometric p-values (HCHG) and other tests, out of $\numberTested$ tested genes from the SCANB data \citep{saal2015sweden}. In all cases, we report on an effect on either side or both sides simultaneously (Algorithm~\ref{alg:1}). 
    Testing for a strictly one-sided effect (Algorithm~\ref{alg:2}) leads to a similar diagram with up to 8 discoveries removed from some groups. The tests Tarone-Ware, Gehan-Wilcoxon, and Feliming-Harrington correspond to different weights in the family of weighted log-rank test \citep{pepe1991weighted}. The family of KONP test is from \cite{gorfine2020k}.
    }
    \label{fig:venn_2sided}
\end{figure}

%% file: table_results_SCANB.tex
\sisetup{
  round-mode   = places, 
  round-precision = 4, 
}

\pgfplotstableset{ 
    highlightrow/.style={
        postproc cell content/.append code={
           \count0=\pgfplotstablerow
            \advance\count0 by1
            \ifnum\count0=#1
            \pgfkeysalso{@cell content/.add={$\bf}{$}}
            \fi
        },
    },
}

\begin{table}
    \begin{center}
    \begin{tikzpicture}        	
    \node[right, scale=.75] (table) at (-1.5,0) {      
    \begin{filecontents*}{csv/results_SCANB.csv}
name,flip,HCHG,Log-rank,Fleming-Harrington01,Tarone-Ware,Peto-Peto
ANKLE2,$<$ med,0.0016,0.168,0.1477,0.3583,0.3445
CLCF1,$<$ med,0.0002,0.2994,0.7082,0.5504,0.5917
DDX5,$<$ med,0.0005,0.2362,0.7622,0.351,0.4129
FAM20B,$<$ med,0.0069,0.436,0.5152,0.9573,0.8929
IL10RB,$>$ med,0.0034,0.0956,0.7959,0.0954,0.1556
MALL,$>$ med,0.004,0.3886,0.3368,0.4743,0.7228
MBD3,$>$ med,0.0074,0.4187,0.1125,0.7458,0.9867
MRAS,$>$ med,0.0005,0.245,0.538,0.262,0.435
MRPS2,$>$ med,0.0024,0.1219,0.1422,0.3103,0.2319
PBX1,$<$ med,0.0,0.2481,0.2918,0.2251,0.4325
\end{filecontents*}
\pgfplotstabletypeset[
      multicolumn names, 
      col sep=comma, 
      display columns/0/.style={
        column name= gene name, 
        column type={l},string type}, 
      display columns/2/.style={
        column name=$\hat{p}(\HCHG)$,
        column type={S[table-format=2.2]},string type},    
    display columns/3/.style={
        column name= $\hat{p}(\LR)$,
        column type={S[table-format=2.2]},string type},
    display columns/4/.style={
         column name= $\hat{p}(\text{FH})$,
         column type={S[table-format=2.2]},string type},
    display columns/5/.style={
        column name= $\hat{p}(\text{TW})$ ,
        column type={S[table-format=2.2]},string type},
    display columns/6/.style={
        column name= $\hat{p}(\text{Peto})$ ,
        column type={S[table-format=2.2]},string type},
    display columns/1/.style={
        column name= increased mortality,
        column type={r},string type},
    every head row/.style={
        before row={\toprule}, 
        after row={
            \midrule} 
            },
    highlightrow/.style={
        postproc cell content/.append code={
            \count0=\pgfplotstablerow
            \advance\count0 by1
            \ifnum\count0=#1
                \pgfkeysalso{@cell content=\textbf{##1}}
            \fi
            },
        },
    every last row/.style={after row=\bottomrule}, 
        ] {csv/results_SCANB.csv}};
\end{tikzpicture}
\caption{
    Several genes with significantly smaller or larger hazards as recognized by higher criticism of hypergeometric P-values (HCHG). At the same time, no effect was recognized by the 
    log-rank (LR) or its variations: Fleming-Harrington (FH), Tarone-Ware (TW), and Peto. The P-value of each test is indicated.  
    The two groups associated with each gene correspond to an expression value higher (lower) than the median. 
    \label{tbl:results_SCANB}}
    \end{center}
\end{table}

%% file: Figs/tkz/fig_real_data.tex
\sisetup{
  round-mode          = places, 
  round-precision     = 3, 
}

\newcommand{\geneNameI}{ANKLE2}
\newcommand{\geneNameII}{CLCF1}
\newcommand{\geneNameIII}{DDX5}

\begin{figure}[ht!]
    \begin{center}
        \begin{tikzpicture}        	
        \node (plot) {\includegraphics[scale = .32]{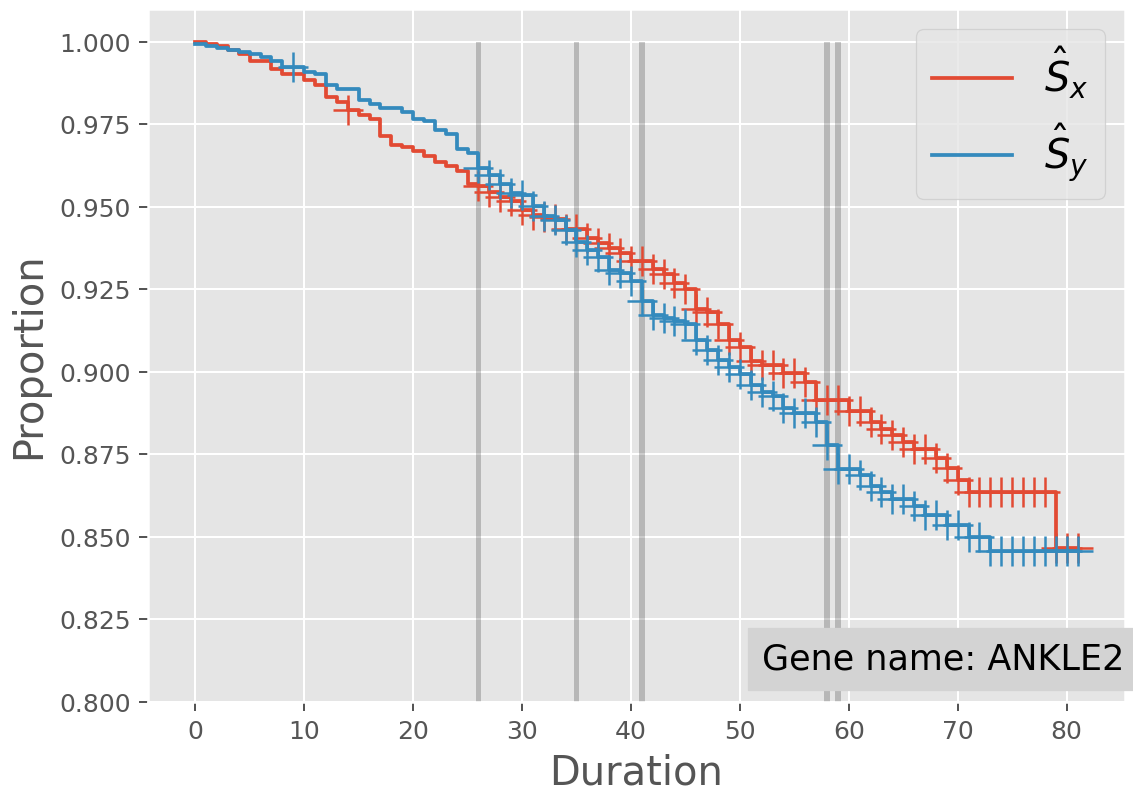}};
        \node[scale=.7, right of=plot, node distance=9cm] {
    \pgfplotstabletypeset[
      multicolumn names, 
      col sep=comma, 
      display columns/0/.style={
        column name= $t$, 
        column type=c,string type},  
      display columns/1/.style={
        column name=$n_x(t-1)$,
        column type=c,string type},
      display columns/2/.style={
        column name=$n_y(t-1)$,
        column type=c,string type},        
    display columns/3/.style={
         column name=  $o_x(t)$,
         column type=c,string type},
    display columns/4/.style={
        column name=  $o_y(t)$,
        column type=c,string type},
    display columns/5/.style={
        column name= $p_t$,
        column type={S},string type},
      every head row/.style={
        before row={\toprule}, 
        after row={
            \midrule} 
            },
        every last row/.style={after row=\bottomrule}, 
        ] {csv/\geneNameI.csv}};
        \end{tikzpicture}
        
        \begin{tikzpicture} 
        \node (plot) {\includegraphics[scale = .32]{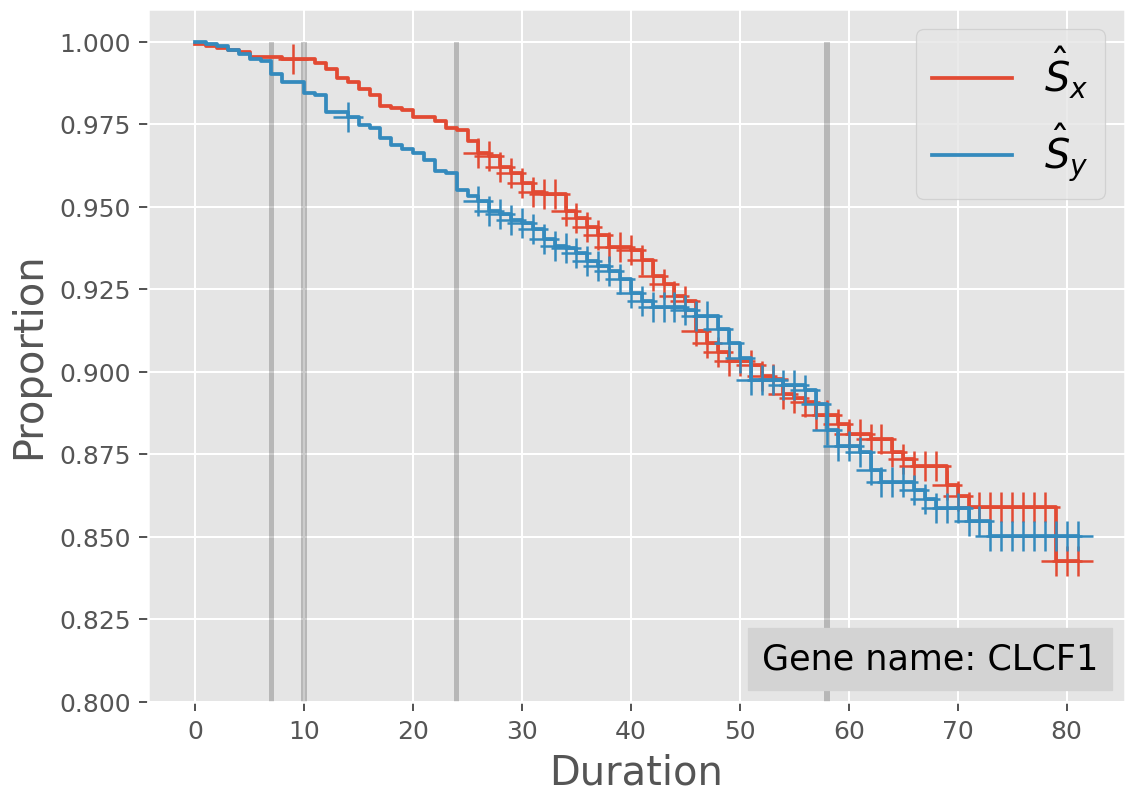}};

        \node[scale=.7, right of=plot, node distance=9cm] {
    \pgfplotstabletypeset[
      multicolumn names, 
      col sep=comma, 
      display columns/0/.style={
        column name= $t$, 
        column type=c,string type},  
      display columns/1/.style={
        column name=$n_x(t-1)$,
        column type=c,string type},
      display columns/2/.style={
        column name=$n_y(t-1)$,
        column type=c,string type},        
    display columns/3/.style={
         column name=  $o_x(t)$,
         column type=c,string type},
    display columns/4/.style={
        column name=  $o_y(t)$,
        column type=c,string type},
    display columns/5/.style={
        column name= $p_t$,
        column type={S},string type},
      every head row/.style={
        before row={\toprule}, 
        after row={
            \midrule} 
            },
        every last row/.style={after row=\bottomrule}, 
        ]{csv/\geneNameII.csv}};
        
        \end{tikzpicture}
\begin{tikzpicture} 
        \node (plot) {\includegraphics[scale = .32]{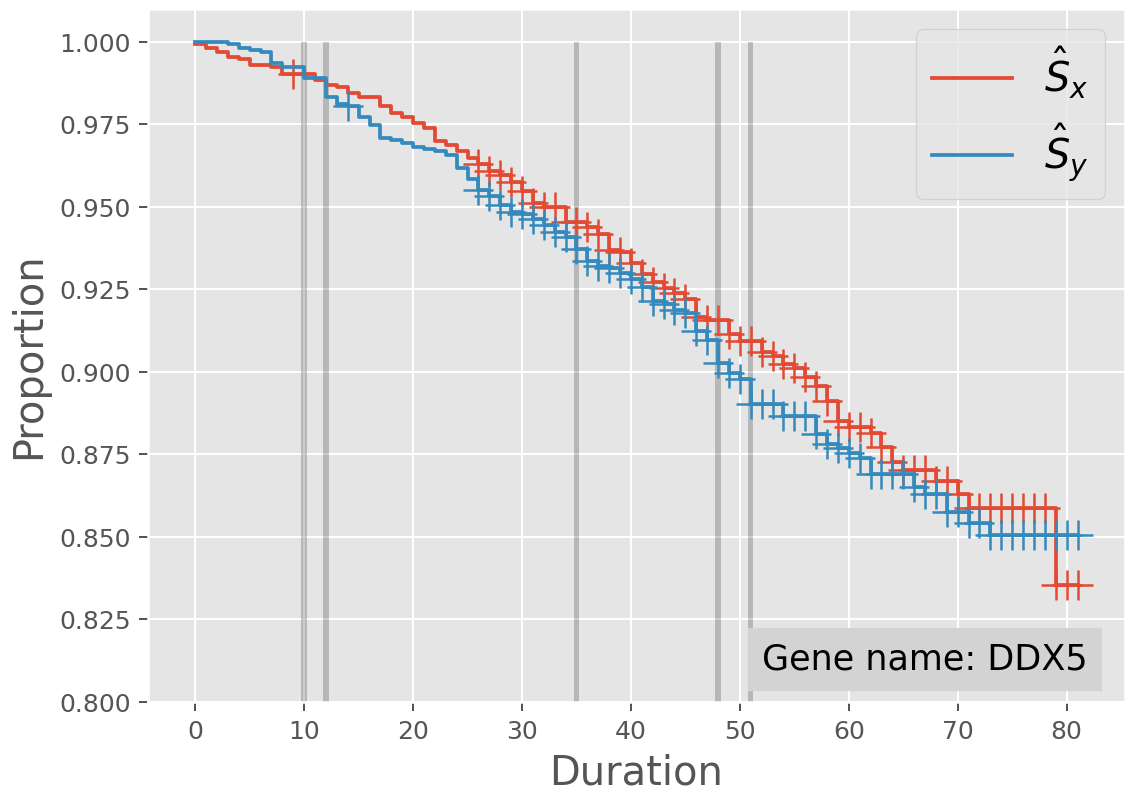}};

        \node[scale=.7, right of=plot, node distance=9cm] {
    \pgfplotstabletypeset[
      multicolumn names, 
      col sep=comma, 
      display columns/0/.style={
        column name= $t$, 
        column type=c,string type},  
      display columns/1/.style={
        column name=$n_x(t-1)$,
        column type=c,string type},
      display columns/2/.style={
        column name=$n_y(t-1)$,
        column type=c,string type},        
    display columns/3/.style={
         column name=  $o_x(t)$,
         column type=c,string type},
    display columns/4/.style={
        column name=  $o_y(t)$,
        column type=c,string type},
    display columns/5/.style={
        column name= $p_t$,
        column type={S},string type},
    every head row/.style={
        before row={\toprule}, 
        after row={
            \midrule} 
            },
        every last row/.style={after row=\bottomrule}, 
        ]{csv/\geneNameIII.csv}};
        \end{tikzpicture}
    \caption{Survival curves in which higher criticism of hypergeometric p-values (HCHG) indicates a significant difference but the log-rank test and some weighted log-rank tests do not. The gray lines correspond to time intervals of suspected excessive hazard, as identified by the set $\Delta^\star$ of \eqref{eq:Delta_def}. The number of at-risk subjects and events in each group at those intervals, and the corresponding hypergeometric P-value, are given in the tables. 
    }
    \label{fig:real_data}
    \end{center}
\end{figure}

%% file: discussion.tex
\subsection{Temporal and time-varying effects}

In previous sections, we discussed the sensitivity of HCHG to detect a global effect 'hiding' in only a few time intervals. Here we emphasize two additional properties of HCHG that are advantageous
in analyzing signals of the temporarily rare hazard departure type compared to other methods. First, as described in Section~\ref{sec:hct}, HCHG offers a mechanism to identify time intervals thought to constitute evidence for a global excessive or reduced hazard. These instances may have important interpretations in some applications. For example, time intervals in which intervention or extra monitoring in medical treatment may be beneficial. Additionally, HCHG can distinguish between effects varying over time that may have opposite trends, e.g., a short-term effect that is different than a long-term one. The studies reported in
\cite{dekker2008survival,johansson2015family,gregson2019nonproportional,daniels2017examining} describe different situations of opposite trends that a procedure based on HCHG can potentially detect, whereas the weighted log-rank statistic \eqref{eq:logrank} might detect only when the effect's time-varying behavior is known before the experiment and thus the weights can be determined correspondingly. 

\subsection{Pooling across time intervals}
When failure events are independent, merging counts over bins of several consecutive time intervals generally reduces the power of a test based on HCHG when the departures are very rare. The intuition here is that the number of non-null instances in one bin is so small that their combined effect diminishes as we increase the bin's size and average the response over its members. This phenomenon is well-understood through previous studies involving rare and weak signal detection models in other settings \citep{arias2011global}. 

When failure events are not independent, pooling across time intervals may improve the detection using HCHG. For example, suppose the presence of an effect causes an increased risk in group $y$ in some period encompassing several consecutive time intervals. In this case, merging counts across bins of time intervals is particularly useful if a bin's size roughly matches the effect duration. When the effect duration is unknown, multi-scale approaches for signal detection might be useful \citep{arias2005near,hall2010innovated,pilliat2023optimal}. In the extreme case when the effects are small and scattered across many instances while the bin size is large, tests based on averaging like the log-rank test or Fisher's combination of the hypergeometric p-values would be preferred.

\subsection{Low risk and rare failures}
The HCHG procedure may be ineffective when the risk of both groups is small such that failure events rarely occur more than once in any given interval. Indeed, this case is very different than the calibration \eqref{eq:calibration_kappa}-\eqref{eq:calibration_rates} in which the number of failure events in each interval goes to infinity. Specifically, the case of few failure events is analogous to the low-counts case of \cite{DonohoKipnis2020} and \cite{arias2015sparse}, which do not correspond to departures on the moderate scale. To our knowledge, testing procedures sensitive to rare departures when the base hazard rate is low is yet an unexplored topic in survival analysis. 

\subsection{Calibration by label permutation}
In Section~\ref{sec:analysis} we used an estimate of the null distribution of $\HCHG$ by label permutation. This calibration is different 
from the theoretical analysis of Section~\ref{sec:analysis} that relied on the asymptotic Brownian bridge behavior of higher criticism \cite{donoho2004higher,shorack2009empirical}. Since the method's power may drop under label permutation calibration, it appears that an asymptotic power characterization of the HCHG under label permutation calibration remains an open challenge. Asymptotic properties of tests based on higher criticism involving permutations have recently been studied in \cite{stoepker2024anomaly,stoepker2023sparse} and in broader contexts in \cite{arias2018distribution,kim2022minimax,dobriban2022consistency}.

Another issue associated with permutation-based calibration stems from a potential mismatch between censoring distributions across groups as discussed in \cite{heimann1998permutational}; see also the discussion in \cite{gorfine2020k}. It is possible that randomizing the hypergeometric tests of \eqref{eq:pvals_def} (e.g., as in \eqref{eq:pval_randomized}) can resolve this issue for calibrating test statistics based on these p-values. We plan to study this point in future work. 


%% file: proofs.tex

\subsection*{Overview}

Theorems~\ref{thm:HC_powerful}, \ref{thm:HC_powerless}, and the region of power reported in Table~\ref{tbl:power}, rely on previous results concerning the ability and impossibility to detect rare mixtures specified in terms of the p-values experiencing moderate non-null departures from \cite{kipnis2021logchisquared}. In the section below, we first state and prove a series of technical lemmas needed to establish the connection between the hypergeometric p-values in our setting \eqref{eq:model_full0}-\eqref{eq:lambda_prime_def} and the rare moderate departures setting of \cite{kipnis2021logchisquared}. The proof of Theorems~\ref{thm:HC_powerful}, \ref{thm:HC_powerless} and \ref{thm:LR_powerless} are provided in subsequent sections.

\subsection*{Asymptotic Notation}

Some technical lemmas concern arrays of real numbers and random variables indexed by $T$ and $t \leq T$ and their asymptotic properties as $T$ goes to infinity. In such cases, we often use the notation $o(1)$ to indicate some deterministic sequence converging to zero uniformly in $t$, and the notation $o_p(1)$ to indicate some sequence of random variables converging to zero in probability uniformly in $t$. We say that a sequence $a(T) \geq b(T)$ eventually if there exists $T_0$ such that $a(T) \geq b(T)$ for all $T\geq T_0$. 

\section{Technical Lemmas}

The following lemma provides an asymptotic lower bound on a hypergeometric P-value of a statistic experiencing moderate deviations. It will be used to establish one side of the convergence in  Lemma~\ref{lem:pvalue_tail_basic} below. 

\begin{lemma}
    \label{lem:pi_lower_bound}
For non-negative integers $x,y,n_x,n_y$ define 
\[
\pi^+(x,y;n_x, n_y) := \Prp{ \HG(n_x + n_y, n_y ,x+y) \geq y + 1 },
\]
where 
\[
\Prp{ \HG(M, N, n) \geq m } = \sum_{k=m}^{n} \frac{\binom{N}{k}\binom{M-N}{n-k}}{\binom{M}{n}}. 
\]
Let $\{\lambda(T)\}$ and $\{a(T)\}$ be positive sequences indexed by $T$ such that, as $T \to \infty$, $a(T) \to \infty$, $a(T) \lambda(T) \geq \log^2(T)$, $a(T)/\log(T) \to 0$, and $\log(\lambda(T))/a(T) \to 0$. For $q>0$ and $T$ sufficiently large, set $\tilde{y}(T,q,x) = \left( \sqrt{x} + \sqrt{q \log(T) - a(T) } \right)^2$. Let sequences $\{n_x(T)\}$ and $\{n_y(T)\}$ obey $\lambda(T)/n_x(T) \to 0$ and $n_x(T)/n_y(T) \to 1$ as $T \to \infty$. There exists $T_0(q)$ such that
\[
\pi^+(x,\tilde{y}(T,q,x);n_x(T), n_y(T)) \geq T^{-q},
\]
for all $T \geq T_0(q)$ and $x \geq \lambda(T) - \sqrt{a(T) \lambda(T)} $. 
\end{lemma}

\subsubsection{Proof of Lemma~\ref{lem:pi_lower_bound}}
\begin{proof}
The proof relies on asymptotic properties of binomial coefficients in the PMF of the hypergeometric distribution. The analysis is similar to \cite[Lemma 5.5]{DonohoKipnis2020}. 
For integers $n$ and $k$ with $n,k \to \infty$ and $k/n \to 0$, Stirling's approximation implies
\begin{align*}
    \binom{n}{k} \sim \left(\frac{n e}{k} \right)^k \frac{1}{\sqrt{2 \pi k}} e^{-\frac{k^2}{2n}(1 + o(1))}. 
\end{align*}
Applying this approximation when $x, y \to \infty$ with $x/n_x = o(1)$ and $y/n_y = o(1)$, we get
\begin{align}
    & \pi^+(x,y;n_x,n_y) = \sum_{k=y+1}^{x+y} \frac{\binom{n_y}{k}\binom{n_x}{x+y-k}}{\binom{n_x + n_y}{x+y}} \geq  \frac{\binom{n_y}{y+1}\binom{n_x}{x-1}}{\binom{n_x + n_y}{x+y}} \nonumber \\
    & \qquad  = \frac{\frac{1}{2\pi \sqrt{(x-1) (y+1)}}\left( \frac{n_y}{y+1}e \right)^{y+1} \exp\left\{-\frac{(y+1)^2}{2n_y}(1+o(1)) \right\} 
      \left( \frac{n_x}{x-1}e \right)^{x-1}  \exp\left\{-\frac{(x-1)^2}{2n_x}(1+o(1)) \right\}} {
     \frac{1}{\sqrt{2\pi (x+y)}}\left( \frac{n_x + n_y}{x+y}e \right)^{x+y} \exp\left\{-\frac{(x+y)^2}{2(n_x+n_y)}(1+o(1)) \right\} 
     } \nonumber  \\
     & = \frac{1}{\sqrt{2\pi}} \sqrt{\frac{x+y}{(x-1)(y+1)}} 
     \left(\frac{n_y^{y+1} n_x^{x-1}}{(n_x+n_y)^{x+y}}\right)
        \left( \frac{x+y}{y+1} \right)^{y+1}
     \left( \frac{x+y}{x-1} \right)^{x-1} \nonumber \\
     & \qquad \times \exp\left\{-\frac{1+o(1)}{2}\left( \frac{(y+1)^2}{n_y} + \frac{(x-1)^2}{n_x} - \frac{(x+y)^2}{n_x+n_y} \right) \right\}. \label{eq:pi_lower_bound:1}
\end{align}
When $n_x/n_y = 1 + o(1)$, we have $n_x + n_y = 2n_y(1+o(1))$ which leads to 
\begin{align}
\frac{n_y^{y+1} n_x^{x-1}}{(n_x+n_y)^{x+y}} = \frac{n_y^{x+y}(1+o(1))^{x-1}}{(2n_y)^{x+y}(1+o(1))^{x+y}} = \left(\frac{1}{2}\right)^{x+y} (1+o(1)).
\label{eq:lemma1:ratio_powers}
\end{align}
For the term in the exponent, assuming $y-x \to \infty$, we get
\begin{align}
    \frac{(y+1)^2}{n_y} + \frac{(x-1)^2}{n_x} - \frac{(x+y)^2}{n_x+n_y} &= \frac{y^2}{n_y} + \frac{x^2}{n_x} - \frac{(x+y)^2}{n_x+n_y} + o\left(\frac{y}{n_y}\right) \nonumber \\
    &= \frac{y^2 n_x(n_x+n_y) + x^2 n_y(n_x+n_y) - (x+y)^2 n_x n_y}{n_x n_y (n_x+n_y)} (1+o(1)) \nonumber \\
    &= \frac{(y n_x - x n_y)^2}{n_x n_y (n_x+n_y)}(1+o(1)) = \frac{(y-x)^2}{2n_y}(1+o(1)).
 \label{eq:pi_lower_bound:2}
\end{align}
Next, consider the term 
\begin{align*}
A & := \left( \frac{x+y}{y+1} \right)^{y+1} \left( \frac{x+y}{x-1} \right)^{x-1} \\
& = -(y+1)\ln\left(1 - \frac{x-1}{x+y}\right) - (x-1)\ln\left(1 - \frac{y+1}{x+y}\right).
\end{align*}
Using the Taylor expansion of $\ln(1-u) = -u - u^2/2 + O(u^3)$, this becomes:
\begin{align}
\label{eq:ln_A}
\ln A & = (x+y)\ln 2 - \frac{(y-x+2)^2}{2(x+y)} + o((y-x)^2/(x+y)^2). 
\end{align}
Combining  \eqref{eq:pi_lower_bound:1}, \eqref{eq:lemma1:ratio_powers}, \eqref{eq:pi_lower_bound:2}, and \eqref{eq:ln_A}, we get
\begin{align}
\log(\pi^+) & \geq -O(\log x) + (x+y)\log(1/2) + \ln A - \frac{(y-x)^2}{4n_y}(1+o(1)) \nonumber \\
& = -O(\log x) - (x+y)\ln 2 + \left( (x+y)\ln 2 - \frac{(y-x+2)^2}{2(x+y)} \right) - \frac{(y-x)^2}{4n_y} + o(1).
 \label{eq:pi_lower_bound:3}
\end{align}
The $(x+y)\ln 2$ terms cancel. With $x^*(T) = \lambda(T) - \sqrt{a(T)\lambda(T)}$, we have
\[
\min_{x \geq \lambda(T) - \sqrt{a(T)\lambda(T)}} \pi^+(x, \tilde{y}(T,x,q);n_x,n_y) =  \pi^+(x^*(T),\tilde{y}(T,x^*(T),q);n_x,n_y). 
\]
We now evaluate the remaining terms for $x=x^*(T)$ and $y = \tilde{y}(T,x^*(T),q)$, under which the conditions $x,y \to \infty$, $y-x \to \infty$, $x/n_x\to 0$, and $y/n_y \to 0$ used previously hold. 

We have:
\begin{align}
\frac{(y-x+2)^2}{2(x+y)} & = \frac{(\tilde{y}-x)^2}{2(x+\tilde{y})}(1+o(1)) \\
& = \frac{4x(q\log T - a(T))}{4x}(1+o(1)) = (q\log T - a(T))(1+o(1)).
\end{align}
Similarly, 
\begin{align}
\frac{(y-x)^2}{4n_y} & = \frac{4x(q\log T - a(T))}{4n_y}(1+o(1)) \\
& = q\log(T) \frac{x}{n_y}(1+o(1)) = o(\log T),
\end{align}
since $x/n_y \sim \lambda(T)/n_y \to 0$. Substituting these into \eqref{eq:pi_lower_bound:3}, we get
\begin{align}
\log(\pi^+) \geq -q\log(T) + a(T) + o(\log(T)) + o(a(T)).
\end{align}
By assumption, $a(T) \to \infty$ faster than $\log(x) \sim \log(\lambda(T))$, thus the dominant terms are $-q\log(T) + a(T)$. Now, consider $T^q \pi^+$:
\[
\log(T^q \pi^+) = q\log T + \log(\pi^+) \geq q\log T + \left(-q\log T + a(T)(1+o(1))\right) = a(T)(1+o(1)).
\]
Since $a(T) \to \infty$, it follows that $\log(T^q \pi^+) \to \infty$, and therefore $T^q \pi^+ \to \infty$. This implies that for any constant $C$, there exists $T_0$ such that for all $T \geq T_0$, $T^q \pi^+ \geq C$. Choosing $C=1$ gives $\pi^+ \geq T^{-q}$, completing the proof. 
\end{proof}

The following lemma provides conditions under which hypergeometric p-values experience moderate departures and characterizes their tail behavior under such departures. 
\begin{lemma}
\label{lem:pvalue_tail_basic}
For non-negative integers $x,y,n_x,n_y$ define 
\[
\pi(x,y;n_x, n_y) = \Prp{ \HG(n_x + n_y, n_y ,x +y) \geq y },
\]
and for $q,s \in \reals$, define $\alpha(q,s) := \left( \sqrt{q} - \sqrt{s} \right)^2$. Let $\bar{\lambda}_t:=\bar{\lambda}(t,T)$,
$n_x := n_x(t,T)$ and $n_y := n_y(t,T)$ be arrays indexed by $T$ and $t \leq T$. Suppose that, as $T \to \infty$, 
$\min_{t\leq T} n_x \to \infty$,  $n_x / n_y = 1 + o(1)$, and $\min_{t \leq T} n_x\bar{\lambda}_t/\log(T) \to \infty$. Define $\bar{\lambda}'_t := \bar{\lambda}'_t(T)=(\sqrt{\bar{\lambda}_t} + \sqrt{\delta}_t)^2$, where $\delta$ satisfies
\[
\delta_t := \delta_t(T) := \frac{r \log(T)}{ 4\frac{n_x n_y}{n_x + n_y} }(1+o(1)). 
\]
Consider the Poisson random variables $X \sim \Pois(\bar{\lambda}_tn_x)$ and $Y \sim \Pois(\bar{\lambda}_t' n_y)$. For any $q > r/2 \geq 0$, 
\begin{align}
\label{eq:lem:pvalue_tail_basic}
     \frac{-\log(\Prp{ \pi(X,Y;n_x, n_y) <  T^{-q})} }{\log(T)} - \alpha(q,r/2)  = o(1).
\end{align}
\end{lemma}

\subsubsection{Proof of Lemma~\ref{lem:pvalue_tail_basic}}
\begin{proof}
The lower bound follows obtained from Lemma~\ref{lem:pi_lower_bound}. 
The upper bound is obtained by a Chernoff-type inequality for the hypergeometric distribution and standard moderate deviation analysis.

For a fixed $x$, $n_x$, $n_y$, and $a>0$. Denote by $y^*(x,a; n_x, n_y)$ the minimal $y$ satisfying
\[
\pi(x,y; n_x, n_y) \leq e^{-a}.
\]
We use the Chernoff inequality for $H \sim \mathrm{HyG}(M,N,n)$ 
\begin{align}
    \label{eq:cheroff}
    \Prp{\sqrt{n} \left(\frac{H}{n} - \frac{N}{M} \right) \geq b } \leq e^{-2 b^2}
\end{align}
in order to bound $y^*(x,a;n_x, n_y)$ (e.g., \eqref{eq:cheroff} follows from \cite[Corollary 1.1]{serfling1974probability}). It follows from \eqref{eq:cheroff} that 
\[
\frac{2}{x+y} \left(y -  (x+y)\frac{n_y}{n_x+n_y}\right)^2 \geq a,
\]
implies $\pi(x,y;n_x, n_y) \leq e^{-a}$. Solving the quadratic expression in this inequality for $y>x\geq 0$, we get
\[
y^*(x, a; n_x, n_y) \geq \frac{\sqrt{8(1-\tilde{\kappa}) a x + a^2 } + 4 \tilde{\kappa}(1-\tilde{\kappa}) x + a }{4 (1-\tilde{\kappa})^2},
\]
where $\tilde{\kappa}:=n_y/(n_x + n_y)$. Setting $a=q \log(T)$, we have
\begin{align}
    & \Prp{\pi(X,Y;n_x, n_y) < T^{-q} }  \geq \Prp{Y \geq y^*(x, a;n_x, n_y)} \nonumber \\
    \quad & = \Prp{ Y \geq \frac{ 4 \tilde{\kappa}(1-\tilde{\kappa}) X + q \log(T) + \sqrt{ 8 (1-\tilde{\kappa})  X q\log(T)  + q^2 \log^2(T)}} {4 (1-\tilde{\kappa})^2}} \nonumber \\
    \quad & = \Prp{ Y \geq \frac{ 4 \tilde{\kappa}(1-\tilde{\kappa}) X \left(1 + \frac{q \log(T)}{X}\right) + \sqrt{ 8 (1-\tilde{\kappa})  X q\log(T) \left(1 + \frac{q \log(T)}{8 (1-\tilde{\kappa}) X}\right) }} {4 (1-\tilde{\kappa})^2}}.
    \label{eq:last_disp}
\end{align}
For any $b \leq \bar{\lambda}_t n_x$,
\begin{align*}
    \Prp{X \leq b} & = \Prp{\bar{\lambda}_t n_x - X \geq \bar{\lambda}_t n_x - b} \\
    & \leq \Prp{\left(\bar{\lambda}_t n_x - X\right)^2  \geq \left(\bar{\lambda}_t n_x - b\right)^2} \\
    & \leq \frac{\bar{\lambda}_t n_x}{\left(\bar{\lambda}_t n_x - b\right)^2};
\end{align*}
the last transition by Markov's inequality applied to the random variable $(\bar{\lambda}_t n_x - X)^2$. Since $\min_{t\leq T} n_x \bar{\lambda}_t / \log(T) \to \infty$, there exists a sequence $b(T)$ such that 
$b(T)/\log(T) \to \infty$, $\max_{t\leq T} b(T)/(\bar{\lambda}_t n_x) \to 0$, and thus 
$\Prp{X \leq b(T)} = o(1)$. Consequently, $X/\log(T)  = o_p(1)$. 
Since $Y$ is independent of $X$, we may replace elements involving $q \log(T)/X$ on the right-hand side of \eqref{eq:last_disp} with the notation $o_p(1)$. By \eqref{eq:last_disp} and $\bar{\kappa}= 1/2 + o(1)$, we obtain
\begin{align}
    & \Prp{\pi(X,Y;n_x, n_y) < T^{-q} }  \nonumber \\ 
    & \quad \geq \Prp{ Y \geq  X \left(1 + o_p(1)\right) + \sqrt{ 4 X q\log(T) \left(1 + o_p(1)\right) } } \nonumber \\
    & \quad = 
    \Prp{ \frac{Y}{2}  \geq  \frac{X}{2} \left(1 + o_p(1)\right) + 2 \sqrt{ \frac{ \frac{X}{2} q\log(T) \left(1 + o_p(1)\right) } {2}}  } \nonumber  \\
    & \quad = 
    \Prp{  \frac{Y}{2} \geq  \left[ \sqrt{ \frac{X}{2}\left(1 + o_p(1)\right)}  +  \sqrt{ \frac{q\log(T) } {2 }} + o_p(1) \right]^2   } \nonumber \\
    & \quad = 
    \Prp{  \sqrt{2 Y } - \sqrt{2 X\left(1 + o_p(1)\right)} + o_p(1) \geq  \sqrt{2 q \log(T)}}.
    \label{eq:pvalue_tail_basic:1}
\end{align}
By the normal approximation to the Poisson and the delta method, the random variables $\sqrt{2X}$ and $\sqrt{2Y}$ are variance stabilized and satisfy
\begin{align*}
   \sqrt{2 X(1+o_p(1))} + o_p(1) & \overset{D}{=} \Ncal(\sqrt{2 \bar{\lambda}_t n_x + o(1) }, 1/2) \overset{D}{=} :\sqrt{2 \bar{\lambda}_t n_x + o(1) } +  Z_x /\sqrt{2} , \\
   \sqrt{2 Y} + o_p(1) & \overset{D}{=} \Ncal(\sqrt{2 \bar{\lambda}'_t n_y }, 1/2 ) \overset{D}{=}: \sqrt{2 \bar{\lambda}'_t n_y } +  Z_y / \sqrt{2}, 
\end{align*}
where $Z_x,Z_y \simiid \Ncal(0,1)$. With $Z\sim \Ncal(0,1)$, we get
\begin{align*}
\sqrt{2 Y} - \sqrt{2 X (1 + o_p(1))} & \overset{D}{=} o_p(1) + Z + 2 \left( \sqrt{ n_y \bar{\lambda}'_t/ 2 } - \sqrt{ n_x \bar{\lambda}_t / 2 }  \right) \\
& = o_p(1) + Z + 2 \sqrt{\frac{n_x n_y}{n_x + n_y} } \left( \sqrt{\bar{\lambda}_t} + \sqrt{\delta_t}  -  \sqrt{\bar{\lambda}_t} \right)  \\
& = o_p(1) + Z + \sqrt{2 \frac{n_x n_y}{n_x + n_y} \cdot 2 \delta_t } = o_p(1) + Z + \sqrt{r \log(T)}.
\end{align*}
From here, a moderate deviation estimate of $\sqrt{2 Y} - \sqrt{2 X}$ as in \eqref{eq:pvalue_tail_basic:1} implies (c.f. \citep{RubinSethuraman1965,zeitouni1998large})
\begin{align}
      \frac{\log \Prp{\pi(X,Y;n_x, n_y) < T^{-q} }}{\log(T)} +o(1) & \geq \lim_{T \to \infty}
    \frac{-\left( \sqrt{2q \log(T)} - \sqrt{r \log(T)} \right)^2}{2 \log(T)} \nonumber \\
    & = -\left(\sqrt{q} - \sqrt{r / 2 } \right)^2 = -\alpha(q,r/2)
    \label{eq:asymp_inequality_1}
\end{align}
whenever $q> r / 2$. 

For the reverse bound, we use
Lemma~\ref{lem:pi_lower_bound} with $\lambda(T) = \bar{\lambda}_tn_x$ and the sequence $a(T) = \log^2(T)/ (\bar{\lambda}_tn_x) + \sqrt{\log(T)}$, which satisfies  $a(T) \to \infty$ as well as the other conditions of Lemma~\ref{lem:pi_lower_bound}. We obtain
\[
2(y^*(x, T^{-q})-1) \leq \left( \sqrt{2 x} + \sqrt{2 q \log(T)(1+o(1))}  \right)^2,
\]
for all $x$ such that $x \geq  n_x \bar{\lambda}_t - \sqrt{a(T) n_x \bar{\lambda}_t}$. Note that $\pi^+(x,y;n_x,n_y) \leq \pi(x,y;n_x,n_y)$ for all $x,y,n_x,n_y$. Therefore, conditioned on the event $A_{t,T} := \{X \geq n_x \bar{\lambda}_t - \sqrt{a(T) n_x \bar{\lambda}_t} \}$, we have
\begin{align*}
\Prp{\pi(X,Y; n_x, n_y) < T^{-q}} & \leq 
\Prp{\pi^+(X,Y; n_x, n_y) < T^{-q}} \\
& \leq \left( \sqrt{2 Y} - \sqrt{2X} + \sqrt{2q \log(T)(1+o_p(1))} \right). 
\end{align*}
From here, the same arguments following \eqref{eq:pvalue_tail_basic:1} above imply 
\[
\frac{\log \Prp{\pi(X,Y;n_x, n_y)< T^{-q} \mid A_{t,T} } }{\log(T)} + o(1) \leq -\alpha(q,r/2). 
\]
Since $a(T) \to \infty$, the normal approximation to the Poisson random variable $X$ and a uniform convergence (Berry Esseen) argument gives
\[
 \Prp{ \frac{X -  n_x \bar{\lambda}_t }{\sqrt{n_x \bar{\lambda}_t}} \geq -\sqrt{a(T)} \mid N_x=n_x } 1 + o(1). 
\]
It follows that $\Prp{A_{t,T}} = 1 + o(1)$, and thus we have the unconditioned asymptotic inequality
\begin{align}
    \label{eq:asump_inequality_2}
\frac{\log \Prp{\pi(X,Y;n_x, n_y)< T^{-q} } }{\log(T)} + o(1) \leq -\alpha(q,r/2). 
\end{align}
Equations \eqref{eq:asymp_inequality_1} and \eqref{eq:asump_inequality_2} imply \eqref{eq:lem:pvalue_tail_basic}. 
\end{proof}
\bigskip

The following lemma estimates the terminal number of subjects in either group under \eqref{eq:model_full1}. We use this lemma to argue that $N_x(T)$ or $N_y(T)$ are not zero infinity often, hence the maximum in \eqref{eq:model_full1} is effectively never `activated'. Later on, we sharpen this result in Lemma~\ref{lem:const_ratio} by showing that $N_x(T)$ and $N_y(t)$ concentrate around the unbounded series
$x_0e^{- \sum_{s \leq T} \bar{\lambda}_s}$ and $y_0e^{-\sum_{s \leq T} \bar{\lambda}_s}$, respectively.

\begin{lemma}
\label{lem:terminal_number_unbounded}
Consider \eqref{eq:model_full0}-\eqref{eq:lambda_prime_def} under the calibration \eqref{eq:calibration_kappa}-\eqref{eq:calibration_rates}. For any $p \geq 1$ and $C>0$, there exists $T_0$ such that  
\begin{align}
    \Prp{N_y(T) \leq C} \leq \Prp{N_x(T) \leq C}  \leq T^{-(p-1)} K^p 
\end{align}
for $T \geq T_0$, where $K$ may depends on $C$ but not on $T$ or $p$. In particular, $\{N_x(T) > 0\}$ and $\{N_y(T) > 0\}$ except for a finite number of $T$ almost surely, by the Borel-Cantelli Lemma.
\end{lemma}
\subsubsection{Proof of Lemma~\ref{lem:terminal_number_unbounded}}
\begin{proof}
First note that $\bar{\lambda}_t' \leq 2 \bar{\lambda}_t$ eventually by \eqref{eq:calibration_delta}, hence it is enough to prove the lemma for $N_x(T)$ obeying \eqref{eq:model_full1} with $\lambda_x = 2 \bar{\lambda}_t$. 

Fix $1 \leq t \leq T$. In what follows, we denote by $\Upsilon_{\lambda}$ an arbitrary random variable with distribution $\Pois(\lambda)$. Given $N_x(t-1)=n$, we have that $N_x(t)=0$ with probability $\Prp{\Upsilon_{n \bar{\lambda}_t} \geq n}$ and $k>0$ with probability $n - \Prp{\Upsilon_{n \bar{\lambda}_t} = n - k}$. By Stirling's approximation, for $k \geq \lambda$, we have $\Prp{\Upsilon_{\lambda} \geq k} \leq e^{-\lambda} (e \lambda/k)^k$. Therefore, 
\begin{align*}
    \Prp{\Upsilon_{\bar{\lambda}k} \geq k} \leq (e \bar{\lambda})^k,
\end{align*}
for any $\bar{\lambda} > 0$, and it follows that for $k-C \geq k\bar{\lambda}_t$, 
\begin{align}
    \Prp{ k - \Upsilon_{k \bar{\lambda}_t } \leq C } = \Prp{ \Upsilon_{k \bar{\lambda}_t } \geq k - C } \leq \left( e\bar{\lambda}_t\right)^{k-C} \leq \left( e M /T\right)^{k-C};
    \label{eq:proof:lemma:terminal:1}
\end{align}
the last transition by \eqref{eq:calibration_rates}. We obtain
\begin{align*}
    \Prp{N_x(t) < C} & = \sum_{k=0}^\infty \Prp{N_x(t) < C \mid N_x(t-1)=k} \Prp{N_x(t-1)=k} \\
    & \leq \sum_{k=0}^{C+p-1} \Prp{N_x(t-1)=k} + \sum_{k=C+p}^\infty \Prp{N_x(t) < C \mid N_x(t-1)=k} \Prp{N_x(t-1)=k} \\
    & \leq \Prp{N_x(t-1) < C + p} + \sum_{k=C+p}^\infty \Prp{ k - \Upsilon_{k \bar{\lambda}t } \leq C } \Prp{N_x(t-1)=k} \\
    & \leq \Prp{N_x(t-1) < C + p} + \sum_{k=C+p}^\infty (eM/T)^{k-C} \Prp{N_x(t-1)=k} \\
    & \leq \Prp{N_x(t-1) < C + p} + \sum_{k=C+p}^\infty (eM/T)^{k-C} \\
    & = \Prp{N_x(t-1) < C + p} + \sum_{j=p}^\infty (eM/T)^{j}
\end{align*}
where the fourth line uses \eqref{eq:proof:lemma:terminal:1} and the final line is a change of index $j=k-C$. For $T>2eM$ we have
\[
\sum_{j=p}^\infty (eM/T)^{j} = \frac{(eM/T)^{p}}{1 - eM/T} \leq \frac{(eM/T)^{p}}{1 - 1/2} \leq 2(eM/T)^{p}.
\]
Therefore, by induction on $t=1,\ldots,T$ for $T > 2eM$, we have
\begin{align*}
\Prp{N_x(T) < C} & \leq 2 T \left(\frac{e M}{T}\right)^p + \one\{{x_0 \leq C + Tp}\} 
\end{align*}
Since \eqref{eq:calibration_rates} implies 
$x_0/(T\log(T)) \to \infty$, and since $2 \leq 2^p$ for $p \geq 1$, the last term is zero and the claim in the lemma follows with $K=(2eM)$. 
\end{proof}

The following lemma says that the size of either group at time $t$ converges in probability to $e^{-\bar{\lambda}t}$, at rate at least $T^{-\beta}$.
\begin{lemma}
\label{lem:const_ratio}
Let $x_0$, $y_0$, $\delta$, $\bar{\lambda}_t$ and $\bar{\lambda}_t'$ be calibrated to $T$ as in  \eqref{eq:calibration_kappa}-\eqref{eq:calibration_rates}. Let $\{N_x(t)\}_{t=1}^T$ and $\{N_y(t)\}_{t=1}^T$ as in \eqref{eq:model_full1}. As $T \to \infty$,
    \begin{align}
        \label{eq:lem:kappa}
     \max_{t\leq T} T \left| \frac{N_x(t)}{x_0 e^{-\sum_{s\leq t} \bar{\lambda}_s }}  - 1 \right| \to 0 \quad \text{and} \quad  \max_{t \leq T} T^{\beta} \left|\frac{N_y(t)}{y_0 e^{ -\sum_{s\leq t} \bar{\lambda}_s}}  - 1\right| \to 0
    \end{align}
    in probability. In particular,
    \begin{align}
    \label{eq:lem:ratio_concentration}
        2\frac{N_x(t) N_y(t)}{N_x(t) + N_y(t)} = x_0 e^{-\sum_{s\leq t} \bar{\lambda}_s}(1 + o(1)) = \frac{x_0+y_0}{2} e^{-\sum_{s\leq t} \bar{\lambda}_s}(1 + o_p(1)),
    \end{align}
    where $o_p(1) \to 0$ uniformly in $t\leq T$ as $T \to \infty$. 
\end{lemma}

\subsubsection{Proof of Lemma~\ref{lem:const_ratio}}
\begin{proof}
We only show the Right-hand side of \eqref{eq:lem:kappa}; to obtain the Left-hand side, replace $\beta$ with $1$ throughout all arguments below. 

Denote 
\[
B_t := N_{y}(t)/(y_0\prod_{s \leq t}(1-\bar{\lambda}_s)). 
\]
Consider the sequence of squared deviations
\[
A_t := \left(B_t-1\right)^2,\qquad t=1,\ldots,T,
\]
and $A_0 = 0$. We will show that this sequence is a submartingale with an expectation that vanishes at the rate advertised in \eqref{eq:lem:kappa}. We first handle the convergence rate of the expectation. For $\lambda>0$, denote $\Upsilon_\lambda \sim \Pois(\lambda)$. Note that for some $b,\lambda,n>0$, we have 
\[
\ex{\left(\frac{n-\Upsilon_{n \lambda}}{b} - 1 \right)^2} = \left( \frac{n(1-\lambda)}{b} - 1 \right)^2 + \frac{n\lambda}{b^2}. 
\]
Given $n=N_y(t-1)>0$, $O_y(t)$ is distributed as $(1-\epsilon)\Pois(n\bar{\lambda}_t) + \epsilon \Pois(n\bar{\lambda}_t')$. By Lemma~\ref{lem:terminal_number_unbounded}, we may assume without loss of generality that given $N_y(t-1) = n$, 
$N_y(t-1)-N_y(t) \sim \Pois(n \bar{\lambda}_t')$ since this holds except for perhaps a finite number of $T$s. We have
\begin{align}
    & \ex{A_t \mid N_y(t-1)} = (1-\epsilon) \left[ A_{t-1} + \frac{\bar{\lambda}_t N_y(t-1) }{y_0^2(1-\bar{\lambda}_t)^{2t}} \right]
    \nonumber \\
    & \quad + \epsilon \left[ \left( \frac{(1-\bar{\lambda}_t')N_y(t-1)}{y_0\prod_{s \leq t}(1-\bar{\lambda}_s)}-1 \right)^2 + \frac{\bar{\lambda}_t' N_y(t-1)}{y_0^2(1-\bar{\lambda}_t)^{2t}} \right] \nonumber \\
    & = A_{t-1} + \left[\bar{\lambda}_t +  \epsilon (\bar{\lambda}_t' - \bar{\lambda}_t) \right]\frac{N_y(t-1) }{y_0^2(1-\bar{\lambda}_t)^{2t}} \nonumber  \\
    & \quad + \epsilon \left[ 
\left( \frac{(1-\bar{\lambda}_t') N_y(t-1)}{y_0\prod_{s \leq t}(1-\bar{\lambda}_s)} - 1 \right)^2 - \left( \frac{(1-\bar{\lambda}_t) N_y(t-1)}{y_0\prod_{s \leq t}(1-\bar{\lambda}_s)} - 1 \right)^2 \right] \nonumber \\
& = A_{t-1} + \left[(1-\epsilon)\bar{\lambda}_t + \epsilon \bar{\lambda}_t' \right]  \frac{B_{t-1}}{y_0(1-\bar{\lambda}_t)^{t+1}} \\
& \qquad + \epsilon \left[\left(\frac{1-\bar{\lambda}_t'}{1-\bar{\lambda}_t} B_{t-1} - 1  \right)^2 - \left(B_{t-1} - 1  \right)^2 \right],
\label{eq:lem_const_ratio:epsilon}
\end{align}
Note that the random variable $B_t$ is bounded from above by $e^{M}$
due to $1 \geq \prod_{s \leq t}(1-\bar{\lambda}_s) \geq e^{-\sum_{s \leq t} \bar{\lambda}_s} \geq e^{-M}$ (here $M$ is from \eqref{eq:calibration_rates}) and $N_y(t-1)\leq y_0$. Likewise, $B_t/\prod_{s \leq t+1}(1-\bar{\lambda}_s) \leq e^{3M}$. Additionally, \eqref{eq:calibration_delta}, \eqref{eq:calibration_rates} imply, 
\begin{align*}
\frac{\delta_t}{\bar{\lambda}_t} & = \frac{\frac{r}{2} \log(T) }{2\bar{\lambda}_t \frac{x_0 + y_0}{2}  e^{-{ \sum_{s \leq t} \bar{\lambda}_s}}} \leq r e^{M} \frac{\log(T) }{2\bar{\lambda}_t (x_0+y_0)} \to 0, 
\end{align*}
hence
\begin{align}
     \eta_t := \frac{\bar{\lambda}_t' - \bar{\lambda}_t}{\bar{\lambda}_t} & = 
     \left\{2 
      \sqrt{\frac{\delta_t}{\bar{\lambda}_t} } + \frac{\delta_t}{\bar{\lambda}_t} \right\}
     = o(1),
     \label{eq:lem_const_ratio:1}
\end{align}
and in particular 
\begin{align}
    \label{eq:lambda_ratio}
    \eta_t = 1 - 
\frac{1 - \bar{\lambda}_t'}{1 - \bar{\lambda}_t} \geq 0.
\end{align}
We obtain
\begin{align*}
    & \ex{A_t \mid N_y(t-1)} = A_{t-1} + o(\bar{\lambda}_t/y_0) + \epsilon o(\bar{\lambda}_t) = A_{t-1} + o(T^{-\beta-1}),
\end{align*}
where in the last transition we used \eqref{eq:calibration_eps} and  \eqref{eq:calibration_rates}. 
Since $A_0 = 0$, by induction on $t=1,...,T$ and \eqref{eq:calibration_rates}, we get
\begin{align}
\label{eq:lem_const_ratio:limit_prop}
    \ex{A_t} = o(T^{-\beta}). 
\end{align}
Next, notice that for $x \to 0$, 
\[
\frac{e^{-x}}{1-x} = 1 + o(x^2).
\]
Since $\sum_{s \leq t} \bar{\lambda}_s^2 \leq M^2/T$ under \eqref{eq:calibration_rates}, we get
\[
1 \leq \frac{e^{- \sum_{s \leq t} \bar{\lambda}_s }}{\prod_{s \leq t}(1-\bar{\lambda}_s)} \leq e^{M^2/T + o(1/T^2)} = 1 + O(1/T).
\]
We conclude that for all $t =0,\ldots,T$,
\[
T^\beta\ex{ \left|\frac{N_y(t)}{y_0 e^{- \sum_{s \leq t} \bar{\lambda}_s}} -1\right|^2} = o(1). 
\]
We now handle the convergence in probability. Denote by $\mathcal F_t$ the sigma-algebra generated by $\{N_y(s),\, s\leq t\}$. We now argue that $\{A_t\}$ is sub-martingale with respect to this filtration. By \eqref{eq:lambda_ratio} and \eqref{eq:lem_const_ratio:1}, we have $\eta_t \geq 2 \sqrt{\delta_{t}/\bar{\lambda}_t}$. By Markov's inequality and \eqref{eq:lem_const_ratio:limit_prop},
\[
\Prp{ \left|1-B_t\right|^2 \geq \eta_t^2/2 } \leq 
\frac{\bar{\lambda}_t}{4\delta_t} \ex{A_t} = O(x_0)o(T^{-\beta}), 
\]
which vanishes due to \eqref{eq:calibration_initial}. This is enough to conclude that the term \eqref{eq:lem_const_ratio:epsilon} is eventually positive and hence $\ex{A_t \mid  \mathcal F_{t-1}} \geq A_{t-1}$. From here, Doob's sub-martingale's inequality (c.f. \cite[P. 870]{shorack2009empirical}) leads to,
\[
\Prp{\max_{t\leq T} A_t \geq T^{-\beta} } \leq T^\beta \ex{A_T} = o(1),
\]
the last transition by \eqref{eq:lem_const_ratio:limit_prop}. All this implies \eqref{eq:lem:kappa}, which also leads to \eqref{eq:lem:ratio_concentration}. 
\end{proof}

The following lemma shows that under the model \eqref{eq:model_full1}, the hypergeometric p-values of \eqref{eq:pvals_def} obey the rare moderate departure  formulation of \cite{kipnis2021logchisquared}.
\begin{lemma}
\label{lem:pvalue_tail}
Let $\bar{\lambda}_t$ and $\delta$ be calibrated to $T$ as in \eqref{eq:calibration_eps}-\eqref{eq:calibration_delta} and $N_x(t)$ and $N_y(t)$ obey \eqref{eq:model_full1}. Suppose that
\[
P_t = p_\HG( \Upsilon'(t); N_x(t) + N_y(t), N_y(t), \Upsilon(t) + \Upsilon'(t) ), 
\]
where, given $N_x(t)$ and $N_y(t)$, $\Upsilon'(t) \sim \Pois(\bar{\lambda}_t'(t)N_y(t))$ and $\Upsilon(t) \sim \Pois(\bar{\lambda}_tN_x(t))$. For $q > r/2 > 0$, we have
\begin{align*}
    \lim_{T \to \infty} \max_{t=1,\ldots T} \left| \frac{-\log(\Prp{ P_t \leq T^{-q}})}{\log(T)} - \alpha(q,r/2) \right| = 0.
\end{align*}
\end{lemma}

\subsubsection{Proof of Lemma~\ref{lem:pvalue_tail}}
\begin{proof}
We show that the conditions of Lemma~\ref{lem:pvalue_tail_basic} hold with probability tending to one as $T\to \infty$.

Define the sequence of events 
\begin{align}
\label{eq:AT_def}
A_T & = \left\{ \min_{t\leq T} N_x(t) \geq x_0 e^{-(M+1)} \right\} \cap  \left\{ \min_{t\leq T}
N_y(t) \geq  y_0 e^{-(M+1)} \right\} \\
& \qquad  \cap 
 \left\{ \max_{t\leq T} \left| \frac{N_x(t)}{N_y(t)} \frac{y_0}{x_0} - 1 \right| <  c_T \right\} \nonumber,
\end{align}
for some positive sequence $c_T$ with $c_T \to 0$ as $T \to \infty$ that will be determined later. Conditioning on $A_T$, Lemma~\ref{lem:pvalue_tail_basic} implies
\begin{align}
\label{eq:lem:pvalue_tail}
    \lim_{T \to \infty} \max_{t=1,\ldots T} \left| \frac{-\log(\Prp{ P_t \leq n^{-q} \mid A_T })}{\log(T)} - \alpha(q,r/2) \right| = 0.
\end{align}
for $q > r/2 > 0$. The claim in the lemma now follows by arguing that $\Prp{A_T} \to 1$. Indeed, by Lemma~\ref{lem:const_ratio},
\begin{align*}
    N_x(t) & = x_0e^{-\bar{\lambda}_tt}\left( 1 + o_p(1) \right) \geq x_0 e^{-M} \left( 1 + o_p(1) \right)\\
    N_y(t) & = y_0 e^{-\bar{\lambda}_tt} \left( 1 + o_p(1) \right) \geq y_0 e^{-M} \left( 1 + o_p(1) \right) 
\end{align*}
hence it follows from \eqref{eq:calibration_rates} that $ \min_{t \leq T} \bar{\lambda}_tN_x(t)/\log(T) \to \infty$ and $ \min_{t \leq T} \bar{\lambda}_tN_y(t)/\log(T) \to \infty$ in probability. Lemma~\ref{lem:const_ratio} also implies that for all $T$ sufficiently large 
\[
  \frac{1-T^{-1}}{1+T^{-\beta}} \leq \frac{N_x(t)}{N_y(t)} \frac{y_0}{x_0} \leq  \frac{1+T^{-1}}{1-T^{-\beta}},\quad t \leq T.
\]
Consequently, 
\[
 -\frac{T^{-\beta}+T^{-1}}{1 + T^{-\beta}} 
 \leq \frac{N_x(T)}{N_y(T)}\frac{y_0}{x_0} - 1  \leq \frac{T^{-\beta}+T^{-1}}{1 - T^{-\beta}} \leq \frac{2T^{-\beta}}{1-T^{-\beta}}
\]
and thus with $c_T := 2T^{-\beta}/(1-T^{-\beta})$ we have 
\[
\Prp{\left\{ \max_{t\leq T} \left| \frac{N_x(t)}{N_y(t)} \frac{y_0}{x_0} - 1 \right| <  c_T \right\}} \to 1.
\]
From here, a union bound on the probability of the complementary event to $A_T$ implies  $\Prp{A_T} \to 1$. 
\end{proof}

The following lemma provides the first two moments of $N_y(t)$ under \eqref{eq:model_full0}-\eqref{eq:lambda_prime_def}. This will be useful in analyzing the asymptotic power of the log-rank test.
\begin{lemma}
\label{lem:logrank_moments}
    Under \eqref{eq:model_full0}-\eqref{eq:lambda_prime_def} and for all $T$ sufficiently large,
    \begin{align}
    \ex{ N_y(t) } & = y_0 \prod_{s=1}^t  
    \left[1 - \left((1-\epsilon) \bar{\lambda}_s + \epsilon \bar{\lambda}_s' \right) \right],
    \label{eq:logrank_mean}
    \end{align}
    and
    \begin{align}
        \Var \left[ N_y(t) \right] & = 
\left[ 
(1-\bar{\lambda}_t)^2 + \epsilon (\bar{\lambda}_t' - \bar{\lambda}_t) \left( \bar{\lambda}'_t + \bar{\lambda}_t - 2 \right)
\right] \Var \left[N_y(t-1)\right] 
\label{eq:logrank_variance} \\
       & \quad +
       \left( \bar{\lambda}_t (1-\epsilon) + \epsilon \bar{\lambda}_t' \right) \ex{N_y(t-1)} \nonumber  \\
       & \quad + \epsilon(1-\epsilon)\left(\bar{\lambda}_t'-\bar{\lambda}_t \right)^2 (\ex{N_y(t-1)})^2.  \nonumber
    \end{align}
\end{lemma}

\subsubsection{Proof of Lemma~\ref{lem:logrank_moments}}
\begin{proof}
By Lemma~\ref{lem:terminal_number_unbounded}, except for perhaps a finite number of $T$'s, we have that
 $N_y(t-1) - N_y(t) = O_y(t) \sim \Pois( \bar{\lambda}'_t N_y(t-1))$ given $N_y(t-1)$. Therefore, the following evaluations of the moments of $N_y(t)$ hold for all $T$ sufficiently large. 
\begin{align*}
\ex{N_y(t) \mid N_y(t-1)} = \left[1 - \left((1-\epsilon) \bar{\lambda}_t + \epsilon \bar{\lambda}_t' \right) \right] N_y(t),
\end{align*}
hence \eqref{eq:logrank_mean} follows by induction on $t$. For \eqref{eq:logrank_variance}, note that
\begin{align*}
    \Var \left[ O_y(t) \mid N_y(t-1) \right] 
    & = \left( \bar{\lambda}_t (1-\epsilon) + \epsilon \bar{\lambda}_t' \right) N_y(t-1) + \epsilon(1-\epsilon)\left(\bar{\lambda}_t'-\bar{\lambda}_t \right)^2N_y^2(t-1)
\end{align*}
where above we used the law of total variance for $O_y(t) \sim (1-\theta)\Pois(\bar{\lambda}_t N_y(t-1)) + \theta \Pois(\bar{\lambda}'_t N_y(t-1))$, $\theta \sim \mathrm{Bernoulli}(\epsilon)$. By the law of total variance, 
\begin{align*}
       \Var \left[ N_y(t) \right] & = 
       \left(1 - (1-\epsilon) \bar{\lambda}_t - \epsilon \bar{\lambda}_t' \right)^2 \Var \left[N_y(t-1)\right] \\
       & +
       \left( \bar{\lambda}_t (1-\epsilon) + \epsilon \bar{\lambda}_t' \right) \ex{N_y(t-1)} + \epsilon(1-\epsilon)\left(\bar{\lambda}_t'-\bar{\lambda}_t \right)^2 \ex{N_y^2(t-1)}.
\end{align*}
Substituting $\ex{N_y^2(t-1)} = (\ex{N_y(t-1) })^2+\Var \left[ N_y(t-1)\right]$ and simplifying leads to \eqref{eq:logrank_variance}. 
\end{proof}

\section{Proof of Theorems}

\subsection{Proof of Theorem~\ref{thm:HC_powerful}
\label{sec:RMD}
}
\begin{proof}
Consider the random variables
\begin{align*}
P_t & = p_{\HG}\left( O_y(t); N_x(t) + N_y(t), N_y(t), O_x(t) + O_y(t) \right), \qquad t=1,\ldots,T,
\end{align*}
where $p_{\HG}$ is defined in \eqref{eq:hyg_pmf}.
Given $N_x(t-1)=n_x(t)$ and $N_y(t-1)=n_y(t)$, $P_t$ is a random variable whose distribution is independent of $P_1,\ldots,P_{t-1}$ and obeys
\begin{align*}
P_t \overset{D}{=} p_{\HG}\left( \Upsilon_y(t); n_x(t) + n_y(t), n_y(t), \Upsilon_x + \Upsilon_y(t) \right), 
\end{align*}
where 
\begin{align*}
\Upsilon_y(t) \sim \Pois(n_y(t) \bar{\lambda}_y(t) ), \qquad \Upsilon_x(t) \sim \Pois(n_x(t) \bar{\lambda}_x(t)), 
\end{align*}
Therefore, considering a sequence of hypothesis testing problems indexed by $T$ and the probability law of $P_t$ given $\{N_x(s), N_y(s)\}_{s \leq t}$, we get the following hypothesis testing problem. 
\begin{align}
\label{eq:hyp_proof}
\begin{split}
    H_0 \,:&\,  P_t \sim \Ucal_t^{(T)}~\text{independently for $t=1,...,T$}, \\
    H_1 \,:&\,  P_t \sim (1-\epsilon)\Ucal_t^{(T)} + \epsilon \Qcal_t^{(T)}~\text{independently for $t=1,...,T$}.
\end{split}
\end{align}
Here $\Ucal_t^{(T)}$ is the distribution of the $t$-th P-value under the null in \eqref{eq:model_full1}, and $\Qcal_t^{(T)}$
is the distribution of 
\[
p_\HG( \Upsilon'(t); n_x(t) + n_y(t), n_y(t), \Upsilon(t) + \Upsilon'(t) ), 
\]
where $\Upsilon'(t) \sim \Pois(\bar{\lambda}'_t n_y(t))$ and $\Upsilon(t) \sim \Pois(\bar{\lambda}_t n_x(t))$. The HC test will turn out to be asymptotically powerful for \eqref{eq:hyp_proof} whenever $r > \rho(\beta)/2$, hence it is also asymptotically powerful for \eqref{eq:model_full1} in this regime. 

Let $U_t \sim \Ucal_t^{(T)}$, $Q_t \sim \Qcal_t^{(T)}$. Since $P_t$ is a P-value under \eqref{eq:model_full1}, the distribution of $P_t$ is super uniform. This is equivalent to
\begin{align}
\label{eq:hyp_proof_U}
     \frac{-\log\Prp{-2 \log(U_t) \geq 2 q \log(T)}}{\log(T)} \leq q, 
\end{align}
for all $t \leq T$ and $T$. In addition, it follows from Lemma~\ref{lem:pvalue_tail} that 
\begin{align}
\label{eq:hyp_proof_Q}
    \lim_{T\to \infty} \max_{t=1,\ldots,T}\left|\frac{-\log\Prp{-2 \log(Q_t) \geq 2 q \log(T)}}{\log(T)} - \alpha(q,r/2) \right| = 0,
\end{align}
with $\alpha(q,s)=(\sqrt{q}-\sqrt{s})^2$.
Equation \eqref{eq:hyp_proof_Q} says that on the moderate deviation asymptotic scale, the sequence $\{-2\log(Q_t)\}_{t=1}^T$ uniformly behaves as a sequence of independent random variables with a noncentral chisquared distribution with one degree of freedom
\[
\chi^2(r/2,1) \overset{D}{=} (Z + \sqrt{r \log(T)})^2,\quad Z \sim \Ncal(0,1).
\]
Hypothesis testing problems involving rare mixtures of
p-values or asymptotic p-values of the form \eqref{eq:hyp_proof} with mixture components
obeying \eqref{eq:hyp_proof_U} and \eqref{eq:hyp_proof_Q} were studied in \cite{kipnis2021logchisquared}. Theorem~\ref{thm:HC_powerful} follows from \cite[Thm. 2]{kipnis2021logchisquared}, and the asymptotic power of tests based on $\mathrm{FDR}^*(p_1,\ldots,p_T)$, $p_{(1)}$, and $F_T$ reported in Table~\ref{tbl:power} follows from 
\cite[Thm. 4-5]{kipnis2021logchisquared}. 
\end{proof}

\subsection{Proof of Theorem~\ref{thm:HC_powerless}}
\begin{proof}
The proof shows that non-null randomized p-values \eqref{eq:pval_randomized} abide by the strong version of the moderate logchisquared approximation for p-values defined in \cite{kipnis2021logchisquared}. The result below, from \cite{kipnis2021logchisquared}, says that all tests based on a rare mixture of such p-values are asymptotically powerless whenever $r < \rho(\beta)$. 
\begin{theorem}{\cite[Cor. 1]{kipnis2021logchisquared}}
\label{thm:strong_logchisq}
Denote\footnote{In the notation of \cite{kipnis2021logchisquared}, $T$ is $n$, $\rho(\beta)$ is $\rho(\beta,1)/2$, and $\chi^2(r)$ is $\chi^2(r/2,1)$.} 
\[
\chi^2(r) \overset{D}{=} \left( Z + \sqrt{ r \log(T)}\right)^2,\qquad Z = \Ncal(0,1).
\]
and denote by $\lognull$ the exponential distribution with mean $2$. Consider the hypothesis testing problem 
\begin{align}
\begin{split}
    H_0^{(T)} & \quad : \quad X_t \sim E_t^{(T)} ,\quad t=1,\ldots,T, \\
    H_1^{(T)} & \quad : \quad X_t \sim (1-T^{-\beta})E_t^{(T)} + T^{-\beta} Q_t^{(T)}, \quad t=1,\ldots,T,
    \label{eq:hyp_log_n_appx}
\end{split}
\end{align}
for some sequences of distributions $Q_t^{(T)}$ and $E_t^{(T)}$. Assume that $Q_t^{(T)}$ is absolutely continuous with respect to $E_t^{(T)}$ for every $t=1,\ldots,T$, for any fixed $q>0$, 
\begin{align}
    \lim_{T\to \infty} \max_{t=1\ldots,T} \frac{\left|\log\left( \frac{dE_t^{(T)}}{d \lognull} (2q\log(T))\right) \right|}{\log(T)} = 0,
    \label{eq:LR_null}
\end{align}
and for any fixed $q > r/2$, 
\begin{align}
    \lim_{T\to \infty} \max_{t=1\ldots,T} \frac{\left|\log\left( \frac{dQ_t^{(T)}}{d \chi^2(r)} (2q\log(T))\right)\right|}{\log(T)} = 0.
    \label{eq:LR_alt}
\end{align}
If $r < \rho(\beta)$, all tests are asymptotically powerless.
\end{theorem}

We use Theorem~\ref{thm:strong_logchisq} with $X_t = -2\log(\tilde{\pi}_t)$ and the hypothesis testing problem \eqref{eq:model_full1}. 
Under the null in \eqref{eq:model_full1}, each randomized P-value $\tilde{\pi}_t$ of \eqref{eq:pval_randomized} has a uniform distribution, hence $-2\log(\tilde{\pi}_t) \sim \lognull$ and \eqref{eq:LR_null} trivially holds. We now show \eqref{eq:LR_alt}.  
We have 
\begin{align*}
\tilde{\pi}(x,y;n_x, n_y) & \geq \Prp{\HG(n_x+n_y, n_y, x+y) > y} \\
& = \Prp{\HG(n_x+n_y, n_y, x+y) \geq y+1} \\
& =: \pi^+(x, y; n_x, n_y),
\end{align*}
hence
\begin{align}
\pi(x;y,n_x,n_y) \geq \tilde{\pi}(x,y;n_x,n_y) \geq \pi^+(x,y;n_x,n_y). 
    \label{eq:converse:proof:sandwitch}
\end{align}
Using the inequality on the Right-Hand of \eqref{eq:converse:proof:sandwitch} and replacing $\pi$ by $\tilde{\pi}$ in the proof of the reverse bound in Lemma~\ref{lem:pvalue_tail_basic}, we obtain 
\begin{align*}
\frac{\log \Prp{\tilde{\pi}(X,Y;n_x, n_y)< T^{-q} } }{\log(T)} + o(1) \leq -
(\sqrt{q} - \sqrt{r/2})^2, \qquad q \geq r/2,
\end{align*}
as $T \to \infty$, 
or
\begin{align}
\Prp{\tilde{\pi}(X,Y;n_x, n_y)< T^{-q} } \leq T^{-(\sqrt{q} - \sqrt{r/2})^2+o(1)},\quad q \geq r/2. 
\label{eq:proof:them:powerless_2}
\end{align}
(notice that \eqref{eq:proof:them:powerless_2} is the counterpart of \eqref{eq:asump_inequality_2} for the randomized p-values.) Replacing $\pi$ with $\tilde{\pi}$ in the proof of Lemma~\ref{lem:pvalue_tail} but otherwise following the exact same steps, we get
\begin{align}
\Prp{\tilde{P}_t < T^{-q}} \leq T^{-(\sqrt{q} - \sqrt{r/2})^2+o(1)},\qquad q \geq r/2,
\label{eq:pval_proof_random_lower}
\end{align}
where $\tilde{P}_t := \tilde{\pi}(\Upsilon'(t); N_x(t)+N_y(t), N_y(t), \Upsilon(t) + \Upsilon'(t))$.
By the Left-Hand side of \eqref{eq:converse:proof:sandwitch} and Lemma~\ref{lem:pvalue_tail}, 
\begin{align}
\Prp{\tilde{P}_t < T^{-q}} \geq T^{-(\sqrt{q} - \sqrt{r/2})^2+o(1)},\qquad q \geq r/2.
\label{eq:pval_proof_random_upper}
\end{align}
Denote by $f_{X_t}(s)$ the density of the random variable $X_t = -2\log(\tilde{p}_t)$. By \eqref{eq:pval_proof_random_lower},  \eqref{eq:pval_proof_random_upper}, and the mean-value theorem, 
\begin{align}
    \label{eq:powerlessness:proof:Q}
\frac{\log\left(f_{X_t}(2q \log(T))\right)}{\log(T)} = -(\sqrt{q}-\sqrt{r/2})^2+o(1).
\end{align}
On the other hand, for $s>0$ we have
\[
\frac{d \chi^2(r)}{ds}(s) = \frac{e^{-\left( \sqrt{s}-\sqrt{r/2}\right)^2}}{2\sqrt{2\pi} \sqrt{s}},  
\]
hence
\begin{align}
    \label{eq:powerlessness:proof:chi}
\frac{\log\left(\frac{d \chi^2(r)}{ds}(2q\log(T) \right)}{\log(T)} = -(\sqrt{q}-\sqrt{r/2})^2+o(1). 
\end{align}
Equations \eqref{eq:powerlessness:proof:Q} and \eqref{eq:powerlessness:proof:chi} implies \eqref{eq:LR_alt}. Theorem~\ref{thm:HC_powerless} follows. 
\end{proof}

\subsection{Proof of Theorem~\ref{thm:LR_powerless}
\label{sec:proof:LR_powerless}}

\begin{proof}
Standard analysis of the log-rank statistics involving independent failure events shows that $\LR$ is asymptotically normal under either hypothesis (c.f. \cite{peto1972asymptotically,schoenfeld1981asymptotic}). It is therefore enough to show that the first two moments of \eqref{eq:logrank} are asymptotically equivalent. 

Under $H_0$, $\LR$ has zero mean and unit variance \citep{peto1972asymptotically}. We evaluate its mean and the variance under $H_1$. Denote 
\[
\bar{\kappa}_t := \frac{N_y(t)}{N_x(t) + N_y(t)}, 
\]
 and notice that
\begin{align*}
    \sum_{t=1}^T O_y(t) - \sum_{t=1}^T E_t & = \sum_{t=1}^T \left( (1-\bar{\kappa}_{t-1})O_y(t) - \bar{\kappa}_{t-1} O_x(t) \right),
\end{align*}
and 
\begin{align*}
    \sum_{t=1}^T V_t & = \sum_{t=1}^T \bar{\kappa}_{t-1}(1-\bar{\kappa}_{t-1})(1+o_p(1)) (O_x(t)+O_y(t))\left( 1 - \frac{O_x(t)+O_y(t)}{N_y(t)}\bar{\kappa}_{t-1} \right),
\end{align*}
Under either hypothesis. Under $H_1$, it follows from Lemma~\ref{lem:const_ratio} that
\begin{align}
\begin{split}
    \label{eq:LR_proof:1}
    N_x(t) & = x_0 e^{-\bar{\lambda}_tt} \left( 1 + o_p(T^{-1}) \right) \\
    N_y(t) & = y_0e^{-\bar{\lambda}_tt}\left( 1 + o_p(T^{-\beta}) \right).
\end{split}
\end{align}
By \eqref{eq:calibration_rates}, \eqref{eq:calibration_kappa}, and \eqref{eq:LR_proof:1},
\[
\bar{\kappa}_t := \frac{1}{2}\left(1 + o_p(T^{-\beta})\right).
\]
Consequently,
\begin{align*}
    \sum_{t=1}^T E_t & = \sum_{t=1}^T \bar{\kappa}_{t-1} (O_x(t)+O_y(t)) = \frac{1+o_p(T^{-\beta})}{2}\sum_{t=1}^T (O_x(t)+O_y(t)) \\
    & = \frac{(1+o_p(T^{-\beta}))}{2} \left[ N_y(0) - N_y(T) +N_x(0)-N_x(T)\right] \\
    & = \frac{(1+o_p(T^{-\beta}))}{2} \left[y_0 \left(1-e^{-\bar{\lambda}_TT}(1+o_p(T^{-\beta}) )\right) \right. \\
    & \qquad \qquad \left. + x_0\left(1-e^{-\bar{\lambda}_TT}(1+o_p(T^{-1}))\right) \right] \\
    & = \frac{y_0 + x_0}{2}(1-e^{-\bar{\lambda}_TT})(1+o_p(T^{-\beta})).
\end{align*}
Similarly, 
\begin{align}
\label{eq:lg_proof_num_mean}
\sum_{t=1}^T O_y(t) = y_0 \left( 1 - e^{-\bar{\lambda}_TT}\right)(1+o_p(T^{-\beta})).
\end{align}
Furthermore, because
\[
\frac{N_x(t-1)}{N_x(t-1)+N_y(t-1)-1}\left( 1 - \frac{O_x(t)+O_y(t)}{N_x(t-1)+N_y(t-1)} \right) = \frac{1 + o_p( T^{-\beta})}{2}(1+o_p(1)),
\]
we have
\begin{align}
\sum_{t=1}^T V_t & = \frac{1+ o_p(1)}{2}\sum_{t=1}^T E_t \nonumber \\
& = \frac{1+ o_p(1)}{2} \frac{y_0 + x_0}{2}(1-e^{-\bar{\lambda}_TT})(1+o_p(T^{-\beta}) ) \nonumber \\
& =  \frac{y_0 + x_0}{4}(1-e^{-\bar{\lambda}_TT})(1+o_p(1))
\label{eq:lg_proof_sum_Vt}
\end{align}
We conclude that
\begin{align*}
\ex{\LR\mid H_1} & =  \sum_{t=1}^T \ex{\frac{ O_y(t)- E_t}
{\sqrt{\sum_{t=1}^T V_t}} \bigm\vert H_1}  \\
& = \frac{ \frac{x_0+y_0}{2} o(T^{-\beta})}{\sqrt{\frac{y_0+x_0}{4}(1-e^{-\bar{\lambda}_T T}) }}(1+o(1)) = o(\sqrt{x_0} T^{-\beta}).
\end{align*}
Additionally, by \eqref{eq:lg_proof_sum_Vt} and similar uses of Lemma~\ref{lem:const_ratio} as above, 
\begin{align*}
    \Var \left[ \LR\mid H_1\right] 
& = \frac{\Var[\sum_{t=1}^T O_y(t) - E_t\mid H_1]}{\frac{y_0+x_0}{4}(1-e^{-\bar{\lambda}_T T})(1+o(1)) } \\
& = \frac{\frac{1}{4}(1+o(1)) \left( \Var \left[\sum_{t=1}^T O_y(t)+O_x(t) \mid H_1\right]  \right)}
{\frac{y_0+x_0}{4}(1-e^{-\bar{\lambda}_T T})(1+o(1))} \\
& = \frac{(1+o(1)) \left(\Var \left[N_x(T)  \mid H_1\right] + \Var \left[ N_y(T) \mid H_1\right] \right) }
{(y_0+x_0)(1-e^{-\bar{\lambda}_T T})(1+o(1))}.
\end{align*}
By Lemma~\ref{lem:logrank_moments} and \eqref{eq:lem_const_ratio:1}, we get 
\[
\Var \left[ N_y(T) \mid H_1 \right] = (1+o(1)) y_0 (1-e^{-\bar{\lambda}_T T})
\]
and likewise for $\Var \left[ N_x(T) \mid H_1 \right] = \Var \left[ N_x(T) \mid H_0 \right]$. It follows that 
\begin{align}
\label{eq:LR:proof2}
\Var[\LR \mid H_0] = 1 =
\Var[\LR \mid H_1](1+o(1)), 
\end{align}
and, since $\ex{\LR \mid H_0} = 0$, 
\begin{align}
    \label{eq:LR:proof1}
\frac{\ex{\LR \mid H_1} - \ex{\LR \mid H_0}}{\sqrt{\Var[\LR \mid H_0]}} = o(\sqrt{x_0} T^{-\beta}). 
\end{align}
By \eqref{eq:calibration_initial}, $\sqrt{x_0} T^{-\beta} = o(1)$ for $\beta>1/2$ hence the first and second moments of $\LR$ are asymptotically equivalent. 
\end{proof}